\input amstex
\documentstyle{amsppt}

\TagsOnRight
\NoBlackBoxes

\loadeusm
\loadeurm

\define\Z{\Bbb Z}
\define\R{\Bbb R}
\define\C{\Bbb C}

\define\Wa{W_{\!+}}
\define\Wm{W_{\!*}}
\define\th{\theta}
\define\T{\Bbb T}

\define\w{\frak{w}}
\define\J{\eusm J}

\define\a{\alpha}
\define\be{\beta}
\define\ga{\gamma}
\define\G{\Gamma}
\define\de{\delta}
\define\la{\lambda}
\define\ep{\epsilon}
\define\e{\varepsilon}
\define\s{\sigma}

\define\U{\Upsilon^{\th,a,b}}
\define\Ut{\Upsilon^\th}
\define\Up{\Upsilon}

\define\dd#1{\dot{#1}}
\define\ddd#1{\ddot{#1}}

\define\tht{\thetag}

\define\I{I}

\define\La{\Lambda}
\define\Lt{\Lambda^\th}

\define\gr{\operatorname{gr}}
\define\const{\operatorname{const}}
\define\Ex{\operatorname{Ex}}
\define\rank{\operatorname{rank}}

\define\ppm{{\prod}^{(+-)}}
\define\ppp{{\prod}^{(++)}}
\define\pp{{\prod}^{(+)}}
\define\pdp{{\prod}^{(+2)}}

\define\eabg{\ep_{\a,\be,\ga}}

\define\pta{\pi_{\tau}}

\define\p{\frak p}
\define\PP{\frak P}

\define\wt{\widetilde}

\rightheadtext{Limits of Jacobi Polynomials}

\topmatter

\title Limits of $BC$--type Orthogonal Polynomials as the
Number of Variables Goes to Infinity
\endtitle

\author
Andrei Okounkov and Grigori Olshanski
\endauthor

\thanks The first author (A.~O.) was partially supported
by NSF and Packard foundation. The second author (G.~O.) was
partially supported by CRDF under grant RM1-2543-MO-03.
\endthanks



\abstract We describe the asymptotic behavior of the multivariate $BC$--type
Jacobi polynomials as the number of variables and the Young diagram indexing
the polynomial go to infinity. In particular, our results describe the
approximation of the spherical functions of the infinite-dimensional symmetric
spaces of type $B,C,D$ or $BC$ by the spherical functions of the corresponding
finite--dimensional symmetric spaces. Similar results for the Jack polynomials
were established in our earlier paper (Intern. Math. Res. Notices 1998, no. 13,
641--682; arXiv: q-alg/9709011). The main results of the present paper were
obtained in 1997.
\endabstract

\address
Department of Mathematics, Princeton University, Princeton, NJ
08544, USA \endaddress \email
okounkov\@math.princeton.edu\endemail
\address Institute for Information Transmission Problems,
Bolshoy Karetny 19, Moscow, 127994, GSP--4, Russia \endaddress
\email olsh\@online.ru \endemail

\toc \head 1. Introduction\endhead

\subhead 1A.  $BC_n$ orthogonal polynomials \endsubhead

\subhead 1B. Statement of the main result \endsubhead

\subhead 1C. Other results \endsubhead

\head 2. Interpolation $BC_n$ polynomials and binomial formula
\endhead

\subhead 2A. Interpolation $BC_n$ polynomials
\endsubhead

\subhead 2B. Binomial formula
\endsubhead

\subhead 2C. Asymptotics of denominators in binomial formula
\endsubhead

\head 3. Sufficient conditions of regularity
\endhead

\head 4. Necessary conditions of regularity
\endhead

\head 5. The convex set $\Ut$ \endhead

\head 6. Spherical functions on infinite--dimensional symmetric
spaces \endhead

\head 7. The $BC_n$ polynomials with $\th=1$ \endhead

\head {} References \endhead

\endtoc
\endtopmatter

\document

\head 1. Introduction
\endhead

In this paper we evaluate the limit behavior of the (suitably
normalized) multivariate orthogonal polynomials associated to
the root system $BC_n$ as $n\to\infty$. For given $n=1,2,\dots$,
these polynomials are viewed as functions on the
$n$--dimensional torus $\T^n$. They are indexed by an arbitrary
partition $\la$ of length at most $n$, and also depend on 3
parameters $\th>0$, $a>-1$, $b>-1$. The normalization is
determined by the condition that the polynomials take value 1 at
the point $(1,\dots,1)$.

We let $n$ go to infinity and assume that the partition
$\la=\la(n)$ varies together with $n$. We obtain necessary and
sufficient conditions on the sequence $\{\la(n)\}$ under which
the corresponding polynomials uniformly converge on any fixed
subtorus $\T^k$, $k=1,2,\dots$. We also describe all possible
limit functions, which live on an infinite--dimensional torus;
it turns out that they depend on countably many continuous
parameters.

The motivation for studying this asymptotic problem comes from
representation theory of infinite--dimensional classical groups.
For certain special values of the parameters $\th,a,b$ (in
particular, $\th$ has to be one of the numbers $\frac12,\,
1,\,2)$ the normalized orthogonal polynomials of type $BC_n$ can
be interpreted as the indecomposable spherical functions on rank
$n$ symmetric spaces of compact type, with restricted root
system $B_n$, $C_n$, $D_n$ or $BC_n$. Likewise, the limit
functions can be interpreted as indecomposable spherical
functions on certain infinite--dimensional analogs of these
symmetric spaces.

Our results not only provide a complete classification of the
spherical functions on infinite--dimensional symmetric spaces
but also explain how these ``infinite--variate'' spherical
functions are approximated by the conventional
(``finite--variate'') spherical functions.

The present paper can be viewed as a continuation of our paper
\cite{OO4} where we studied a similar asymptotic problem for the
Jack polynomials. The Jack polynomials  depend on a single
continuous parameter $\th>0$. For 3 special values $\frac12,
1,\,2$ of this parameter, the (suitably normalized) Jack
polynomials with $n$ variables can be interpreted as
indecomposable spherical functions on rank $n$ compact symmetric
spaces with restricted root system $A_n$:
$$
U(n)/O(n), \qquad (U(n)\times U(n))/U(n), \qquad U(2n)/Sp(n).
$$

In the case $\th=1$ the Jack polynomials become the Schur
polynomials. The normalized Schur polynomials can also be viewed
as the normalized irreducible characters of the groups $U(n)$.
The large $n$ asymptotics of these characters was found, for the
first time, in the pioneer work of A.~M.~Vershik and
S.~V.~Kerov, \cite{VK}.

The results of our paper \cite{OO4} and of the present paper
together provide a far generalization of \cite{VK}, involving
all families of multivariate orthogonal polynomials connected
with classical root systems. A natural field of application for
our results is infinite--dimensional noncommutative harmonic
analysis in the spirit of the papers Olshanski \cite{O4} and
Borodin--Olshanski \cite{BorO}. Those papers present a detailed
study of harmonic analysis on the infinite--dimensional unitary
group, and they substantially use the large $n$ asymptotics of
Schur polynomials. One can expect that a similar theory can be
built for other infinite--dimensional classical groups or
symmetric spaces, which will necessarily imply a similar use of
more general orthogonal polynomials.

We proceed to a more detailed description of the results
obtained in the present paper.

\subhead 1A. $BC_n$ orthogonal polynomials
\endsubhead

Throughout the paper $n=1,2,\dots$ denotes a natural number. Let
$W$ denote the $BC_n$ Weyl group
$$
W=S(n)\ltimes \Z_2^n\,.
$$
We shall need two copies $\Wm$ and $\Wa$ of this group acting on
functions in $n$ variables. The $S(n)$ part in both cases
permutes the variables $z_1,\dots,z_n$,  and the $\Z_2^n$ part
acts by
$$
f(z)\mapsto f(z^{\pm1}_1,\dots,z^{\pm1}_n)\,, \quad f(z)\mapsto
f(\pm z_1,\dots,\pm z_n)\,,
$$
in $\Wm$ and $\Wa$, respectively.

Given a partition $\la$, we denote by $\ell(\la)$ the number of
nonzero parts of $\la$. The $BC_n$ orthogonal polynomials (also
called the $BC_n$ Jacobi polynomials) are certain
$\Wm$--invariant Laurent polynomials in $n$ variables, labelled
by partitions $\la$ with $\ell(\la)\le n$ and depending on 3
parameters
$$
\th> 0\,, \quad a,b> -1\,.
$$
The polynomials are defined as follows. The parameters $\th,a,b$
specify an inner product,
$$
(f,g) = \int_{\T^n} f(z)\, \overline{g(z)} \,\w(z)
\cdot\operatorname{Haar}(dz)\,,
$$
of functions on the
$n$-dimensional torus
$$
\T^n=\{(z_1,\dots,z_n)\subset\C^n\}\,, \quad |z_i|=1\,,
$$
where ``$\operatorname{Haar}$'' is the Haar measure on $\T^n$
and $\w(z)$ is the following $\Wm$-invariant weight function
$$
\w(z)=\prod_{1\le i<j \le n} |z_i-z_j|^{2\th} |1-z_i z_j|^{2\th}
\prod_{1\le i\le n} |1-z_i|^{2a+1} |1+z_i|^{2b+1} \,. \tag1.1
$$

The polynomials in question, denoted as $\J_\la(z;\th,a,b)$, are
$\Wm$-invariant, orthogonal with respect to the above inner
product, and satisfy the triangularity condition
$$
\J_\la(z;\th,a,b)=z^\la+\dots\,,  \tag1.2
$$
where $z^\la=z_1^{\la_1}\cdots z_n^{\la_n}$ and dots stand for
lower monomials in the lexicographic order. These  properties
characterize the polynomials $\J_\la(z;\th,a,b)$ uniquely.

Actually, the polynomials $\J_\la(z;\th,a,b)$  possess a
stronger triangularity property. To state it we need some
notation. Given a partition $\mu$ with $\ell(\mu)\le n$, set
$$
\wt m_\mu(z)=\sum_{\nu\in W_+(\mu)}z^\nu, \qquad z\in\T^n
\tag1.3
$$
(summed over weights $\nu$ in the $W_+$--orbit of $\mu$). This
is an analog of the monomial symmetric function for the root
system $BC_n$. Next, let $\e_1,\dots,\e_n$ be the canonical
basis in $\Z^n$. Write $\mu\ll\la$ if the vector
$(\la_1-\mu_1,\dots,\la_n-\mu_n)\in\Z^n$ can be written as a
linear combination of the vectors $\e_i-\e_j$ ($1\le i<j\le n$)
and $\e_i$ ($1\le i\le n$) with nonnegative integral
coefficients. In this notation, the refinement of the
triangularity condition  has the form:
$$
\J_\la(z_1,\dots,z_n;\th,a,b)
=\sum_{\mu\ll\la}u_{\la\mu}(\th,a,b)\,\wt m_\mu(z_1,\dots,z_n),
\tag1.4
$$
where $u_{\la\mu}(\th,a,b)$ are certain coefficients such that
$u_{\la\la}(\th,a,b)=1$.

This  is a specialization of the general definition of the
multivariate Jacobi polynomials corresponding to an arbitrary
root system $R$. See, for example, Heckman's lectures published
in \cite{HS},  Macdonald \cite{Ma4}, and Koornwinder's
expository paper \cite{K2}. Those general polynomials depend on
a Weyl group invariant function
$$
\a\mapsto k_\a\,, \quad \a\in R\,,
$$
on the root system $R$. In our case
$$
R=\left\{\pm\e_i\pm\e_j\,,\pm\e_i\,,\pm2\e_i\right\}\,,
$$
where $\{\e_i\}$, as above, is the standard basis of $\R^n$, and
$$
k_{\pm\e_i\pm\e_j}=2\th\,, \quad k_{\pm\e_i}=a-b\,, \quad \quad
k_{\pm2\e_i}=b+\tfrac12\,. \tag1.5
$$
The numbers $k_\a$ are viewed as formal root multiplicities.
Thus, the ``half--sum of the positive roots'' is defined as
$$
\rho:=\frac12\sum_{\a>0} k_\a \a =
(\th(n-1)+\s,\,\th(n-2)+\s,\,\dots,\,\th+\s,\,\s)\,,
$$
where
$$
\s=\frac{a+b+1}2\,.
$$

We shall need the following known fact:

\proclaim{Proposition 1.1} Assume $a\ge b\ge-\tfrac12$ so that
the formal multiplicities defined in \tht{1.5} are nonnegative.
Then in the expansion \tht{1.4}, the coefficients
$u_{\la\mu}(\th,a,b)$ are all nonnegative.
\endproclaim

\demo{Proof} See Macdonald's paper \cite{Ma4}, formula
\tht{11.15} and the argument following it. \qed
\enddemo

By an appropriate specialization of the parameters $a$, $b$ one
can obtain the orthogonal polynomials associated to the root
systems $B_n$, $C_n$ or $D_n$. Namely, we have to set
$b=-\tfrac12$, $a=b\ne-\tfrac12$ or $a=b=-\tfrac12$,
respectively. It should be noted, however, that in the $D_n$
case, the polynomial $\J_\la(z;\th,
a,b)=\J_\la(z;\th,-\tfrac12,-\tfrac12)$ with $\la_n>0$ is the
sum of certain two ``twin'' $D_n$--polynomials.

Note that the orthogonal polynomials for the $A_n$ root system
are very closely related to Jack polynomials (see Beerends and
Opdam \cite{BeO}); the asymptotics of Jack polynomials as
$n\to\infty$ was studied in our paper \cite{OO4}. As for the
exceptional root systems, there are, obviously, no $n\to\infty$
asymptotic problems.

Similarly to the Jack polynomials, the polynomials
$\J_\la(z;\th,a,b)$ are eigenfunctions of $n$ commuting
differential operators  and describe excitations in certain
completely integrable quantum many body systems which were
introduced by Olshanetsky and Perelomov \cite{OP}.

For certain special values of the parameters $\th,a,b$, the root
data $(R,\{k_\a\})$ given by \tht{1.5} correspond to the
restricted root system of a rank $n$ compact symmetric space
$G(n)/K(n)$ of classical type (for more detail, see \S6). Then
the commuting differential operators mentioned above become the
radial parts of invariant differential operators on $G(n)/K(n)$.

Next, we introduce the {\it normalized} Jacobi polynomials
$$
\Phi_\la(z;\th,a,b):=\frac{ \J_\la(z_1,\dots,z_n;\th,a,b)}
{\J_\la(\,\underbrace{1,\dots,1}_n\,;\th,a,b)}\,. \tag1.6
$$
For special values of $\th,a,b$ these are the {\it
indecomposable spherical functions\/} on the corresponding
symmetric space $G(n)/K(n)$.

The Jacobi polynomials $\J_\la(z;\th,a,b)$ are degenerations as
$q\to 1$ of the 6--pa\-ra\-me\-tric Koornwinder polynomials
which are eigenfunctions of certain commuting $q$-difference
operators (see Koornwinder \cite{K1}, van Diejen \cite{Di}).

Finally, note that the Jacobi polynomials $\J_\la(z;\th,a,b)$
can be transferred from the torus $\T^n$ to the cube $[-1,1]^n$
via the map
$$
z_i \mapsto x_i=\frac{z_i+z^{-1}_i}2\,,
$$
which takes $\Wm$-invariant polynomials in $z$ to symmetric
polynomials in $x$. The weight function \tht{1.1} is then
replaced by the weight function
$$
\prod_{1\le i<j\le n} |x_i-x_j|^{2\th} \prod_{1\le i \le n}
(1-x_i)^a (1+x_i)^b
$$
with respect to Lebesgue measure $dx$ on the cube. When $n=1$,
this is the familiar weight function for the classical Jacobi
polynomials.

\subhead 1B. Statement of the main result
\endsubhead

Fix some $\th>0$. We could have also fixed some $a,b>-1$;
however, at no extra cost, we can consider the following more
general situation. Namely, we fix two sequences
$$
\{a_n\}\,, \{b_n\}\,,
$$
such that $a_n , b_n > -1$ and the limits
$$
\bar a=\lim_{n\to\infty} \frac{a_n}n\,, \qquad \bar
b=\lim_{n\to\infty} \frac{b_n}n \tag1.7
$$
exist. In particular, if $a_n$ and $b_n$ do not depend on $n$
then
$$
\bar a=\bar b=0.
$$
We consider $\th$, $\{a_n\}$, and  $\{b_n\}$ as fixed parameters
of our problem and study the limit behavior as $n\to\infty$ of
the functions
$$
\Phi_{\la(n)} (z_1,\dots,z_k,\underbrace{1,\dots,1
}_{\text{$n-k$ times}};\th,a_n,b_n)\,,   \tag1.8
$$
where $k=1,2,3,\dots$ is fixed and  $\{\la(n)\}$ is a sequence
of partitions with $\ell(\la(n))\le n$,
$$
\la(n)=(\la(n)_1\ge\la(n)_2\ge \dots \ge \la(n)_n\ge 0). \tag1.9
$$

\definition{Definition 1.2} Let $\{\la(n)\}$ be a sequence of partitions
as in \tht{1.9}.

(i) $\{\la(n)\}$ is said to be {\it regular\/} if for every
fixed $k$, the functions \tht{1.8} uniformly converge on the
torus $\T^k$, as $n\to\infty$.

(ii) $\{\la(n)\}$ is said to be {\it infinitesimally regular\/}
if  for every fixed $k$, the Taylor expansions of \tht{1.8} in a
local system system of coordinates about the point
$(1,\dots,1)\in\T^k$  have a coefficient--wise limit; in this
case we shall say that the functions \tht{1.8} {\it converge
infinitesimally\/}. In other words, this type of convergence
means the convergence of jets at the point $(1,\dots,1)\in\T^k$.

(iii) $\{\la(n)\}$ is said to be {\it minimally regular\/} if
the functions
$$
\Phi_{\la(n)}(z,\underbrace{1,\dots,1 }_{\text{$n-1$
times}};\th,a_n,b_n)\,, \quad |z|=1\,,
$$
converge pointwise to a continuous function on $\T^1$.
\enddefinition

Given a partition $\la$, let $\la'$ denote the conjugate
partition. That is, $\la'_i$ is the length of the $i$th column
in the diagram of $\la$. Let $|\la|$ denote the sum of the parts
of $\la$ (equivalently, the number of the boxes in the
corresponding diagram).

\definition{Definition 1.3} A sequence \tht{1.9} is said to be
a {\it Vershik-Kerov sequence\/} (VK sequence, for short)  if
the following limits exist:
$$
\aligned
\a_i &:= \lim \frac{\la(n)_i}n \quad<\infty\,,\quad i=1,2,\dots\,,\\
\be_i &:= \lim \frac{(\la(n))'_i}n \quad<\infty\,,\quad
i=1,2,\dots\,, \\
\de &:= \lim \frac{|\la(n)|}n \quad<\infty\,.
\endaligned \tag1.10
$$
It is readily checked (see below) that the number
$$
\ga:=\de - \sum (\a_i + \be_i)
$$
is nonnegative. The numbers $\a_i$, $\be_i$, $\ga$ (or $\a_i$,
$\be_i$, $\de$) are called the {\it VK parameters\/} of the
sequence $\{\la(n)\}$.
\enddefinition

Let us give a slightly different (but equivalent) definition of
the VK parameters, which makes evident the inequality $\ga\ge0$.
Let $d(n)$ denote the number of diagonal boxes in the Young
diagram corresponding to $\la(n)$. For $i=1,\dots,d(n)$, we
replace in \tht{1.10} the row lengths $\la(n)_i$ and the column
lengths $(\la(n))'_i$ by the respective {\it modified Frobenius
coordinates\/}
$$
\la(n)_i-i+\tfrac12,\qquad (\la(n))'_i-i+\tfrac12,
$$
and for $i>d(n)$, we replace the row and column lengths by
zeros, which does not affect the definition of $\a_i$ and
$\be_i$. On the other hand, the sum of the modified Frobenius
coordinates equals $|\la(n)|$. After the limit transition, this
turns into the inequality $\sum(\a_i+\be_i)\le\de$, so that
$\ga\ge0$.

The following is the main result of the present paper.

\proclaim{Theorem 1.4} Let $\th>0$ and $a_n\ge b_n\ge -\tfrac12$
be parameters satisfying \tht{1.7}. Let
$\{\la(n)\}_{n=1,2,\dots}$ be a sequence of partitions with
$\ell(\la(n))\le n$.

{\rm(i)} All 3 regularity properties of Definition 1.2 are
equivalent to each other and  are also equivalent to the VK
conditions \tht{1.10} of Definition 1.3.

{\rm(ii)} If $\{\la(n)\}$ is a VK sequence with parameters
$\a_i$, $\be_i$, $\ga$ then for any fixed $k=1,2,\dots$
$$
\lim_{n\to\infty} \Phi_{\la(n)}
(z_1,\dots,z_k,\underbrace{1,\dots,1 }_{\text{$n-k$
times}};\th,a_n,b_n)=\prod_{j=1}^k \phi_{\a,\be,\ga,\bar a,\bar
b}\left(\frac{z_j+z^{-1}_j}2\right)\,,
$$
where $\bar a, \bar b$ are defined in \tht{1.7}, and
$\phi_{\a,\be,\ga,\bar a,\bar b}$ is the following function of a
single variable $x\in[-1,1]$
$$
\phi_{\a,\be,\ga,\bar a,\bar b}(x)=e^{\ga(x-1)}
\prod_{i=1}^{\infty} \dfrac{
1+\dfrac{\be_i}2\left(\dfrac{2\th+\bar a+\bar
b-\th\be_i}{\th+\bar a}\right) (x-1)}
{\left(1-\dfrac{\a_i}{2\th}\left(\dfrac{2\th+\bar a+\bar
b+\a_i}{\th+\bar a}\right) (x-1)\right)^\th}\,. \tag1.11
$$
In particular, if $\bar a=\bar b=0$ (which is the case when
$a_n$ and $b_n$ do not depend on $n$) then the above expression
can be written as
$$
\multline \phi_{\a,\be,\ga}\left(\frac{z+z^{-1}}2\right)\\
=e^{\frac\gamma2(z+z^{-1}-2)} \prod_{i=1}^\infty
\frac{\left(1+\frac{\beta_i}2(z-1)\right)
\left(1+\frac{\beta_i}2(z^{-1}-1)\right)}
{\left(\left(1-\frac{\alpha_i}{2\theta}(z-1)\right)
\left(1-\frac{\alpha_i}{2\theta}(z^{-1}-1)\right)\right)^\theta}\,.
\endmultline
\tag1.12
$$
\endproclaim

Note that the infinite products in these formulas are convergent
because $\sum(\a_i+\be_i)<\infty$.

Note also that there exists an {\it a priori\/} argument (see
\cite{O2, \S23} and \cite{O3}) explaining why the limit
functions in Theorem 1.4 factorize. This argument works for
special values of $(\th,a,b)$ when the limit functions admit a
representation theoretic interpretation, see \S1C below.

The proof of Theorem 1.4 is completed in section 4. The strategy
of the proof is the same as in our paper \cite{OO4}. Our main
technical tool is the {\it binomial formula\/} for the Jacobi
polynomials $\J_\la(z;\th,a,b)$ which involves the so called
interpolation $BC_n$ polynomials $\I_\mu(x)$. These objects are
discussed in detail in section 2 below. Note that the binomial
formula we need is a degeneration of a more general binomial
formula for the Koornwinder polynomials obtained in \cite{Ok1}.

There is, however, a difference in the way we establish the
sufficient conditions of regularity in the Jack and Jacobi
cases. The argument in the Jack case relies heavily on the use
of generating series for one--row shifted Jack polynomials. It
is unknown how to evaluate such a series in the Jacobi case.
Instead, there is a simple argument available which uses the
stability of the polynomials $\I_\mu$.

\subhead 1C. Other results \endsubhead

In section 5, we fix parameters $a,b$ such that $a\ge
b\ge-\frac12$. For $n=1,2,\dots$ let $\U_n$ be the convex set of
functions on the torus $\T^n$ which can be written as convex
linear combinations of functions
$\Phi_\la(z_1,\dots,z_n;\th,a,b)$ with $\ell(\la)\le n$. It
turns out that the specialization $z_n=1$ determines an affine
map $\U_n\to\U_{n-1}$. Consequently we can form the projective
limit of the convex sets $\U_n$ as $n\to\infty$, which we denote
as $\Ut$. Consider the set $\T^\infty_0=\varinjlim\T^n$ whose
elements are infinite vectors
$(z_1,z_2,\dots)\in\T\times\T\times\dots$ with finitely many
coordinates $z_i$ distinct from 1. Then elements of $\Ut$ can be
described as functions $\varphi(z_1,z_2,\dots)$ on $\T^\infty_0$
such that for any $n=1,2,\dots$, the function
$$
\varphi_n(z_1,\dots,z_n)=\varphi(z_1,\dots,z_n,1,1,\dots)
$$
on $\T^n$ belongs to $\U_n$.

It is clear that $\Ut$ is a convex set. In Theorem 5.2 we show
that there is a one--to--one correspondence between extreme
points of $\Ut$ and collections $(\a,\be,\ga)$ of VK parameters,
$$
\a=(\a_1\ge\a_2\ge\dots\ge0), \quad
\be=(\be_1\ge\be_2\ge\dots\ge0), \quad \ga\ge0, \quad
\sum(\a_i+\be_i)<\infty.
$$
Given $(\a,\be,\ga)$, the corresponding extremal function on
$\T^\infty_0$ has the form
$$
\varphi(z_1,z_2,\dots) =\prod_{j=1}^\infty \phi_{\a,\be,\ga}
\left(\frac{z_j+z_j^{-1}}2\right),
$$
where the functions $\phi_{\a,\be,\ga}$ are those defined in
\tht{1.12}. That is, the possible limits of normalized Jacobi
polynomials are precisely the extreme points of $\Ut$. In
particular, this implies that the set $\Ut$ does not depend on
parameters $a,b$.

In section 6, we consider 7 infinite--dimensional symmetric
spaces $G/K$ defined as the inductive limits of classical
compact symmetric spaces $G(n)/K(n)$ of type $B$, $C$, $D$ or
$BC$. We assume that parameters $\th,a,b$ take special values
depending on the series $\{G(n)/K(n)\}$. Then the convex set
$\Ut$ can be identified with the set of positive definite,
two--sided $K$--invariant, normalized functions on $G$. As a
corollary of Theorem 5.2 we obtain an explicit description of
the indecomposable spherical functions on $G/K$. Moreover,
Theorem 1.4 shows how these functions are approximated by the
indecomposable spherical functions of the rank $n$ compact
symmetric spaces $G(n)/K(n)$.

An important particular case of this result, corresponding to
the large $n$ asymptotic behavior of the characters of the
orthogonal and symplectic groups,  was earlier obtained by Boyer
\cite{Boy}.

In section 7 we present an elementary derivation of two basic
facts about the Jacobi polynomials $\J_\la$, the binomial
formula and the branching rule, for the case $\th=1$.

\remark{Acknowledgment} We are grateful to the referee for
valuable remarks.
\endremark

\head 2. Interpolation $BC_n$ polynomials and binomial formula
\endhead

\subhead 2A. Interpolation $BC_n$ polynomials
\endsubhead

Fix $\th>0$ and denote by $\Lt_n$ the  algebra of polynomials in
$n$ variables $x_1,\dots,x_n$, symmetric in variables $x_i-\th
i$, $i=1,\dots,n$. Consider the projective limit of these
algebras
$$
\Lt=\varprojlim \Lt_n
$$
taken in the category of filtered algebras, with respect to
homomorphisms sending the last variable to 0. Here the
filtration is defined by the total degree $\deg(\,\cdot\,)$ of a
polynomial. That is, an element $f\in\Lt$ is a sequence of
polynomials $f_n\in\Lt_n$ such that
$f_{n+1}(x_1,\dots,x_n,0)=f_n(x_1,\dots,x_n)$ and $\deg (f_n)$
remains bounded. Elements $f\in\Lt$ can be evaluated at any
infinite vector $x=(x_1,x_2,\dots)$ with finitely many nonzero
coordinates. In particular, for any partition $\la$, the value
$f(\la)$ is well defined. Let $h$ be one more variable and
consider the algebras
$$
\Lt(h)=\Lt\otimes\C(h)\,, \quad \Lt[h]=\Lt\otimes\C[h]\,,
$$
and define $\Lt_n(h)$ and $\Lt_n[h]$ similarly.

\proclaim{Proposition 2.1} {\rm(i)} In the $\C(h)$--algebra
$\Lt_n(h)$ defined above, there exist polynomials
$\I_\mu(x_1,\dots,x_n;\th;h)$ indexed by arbitrary partitions
$\mu$ with $\ell(\mu)\le n$, satisfying the following Newton
interpolation conditions:

{\baselineskip=15 pt \roster \item
$\deg(\I_\mu(\,\cdot\,;\th;h))=2|\mu|$, \item
$\I_\mu(x_1,\dots,x_n;\th;h)$ is $\Wa$--invariant in variables
$x_i-\th i+ h$, $i=1,\dots,n$, \item $\I_\mu(\la;\th;h)=0$ if
$\la=(\la_1,\dots,\la_n)$ is a partition with $\ell(\la)\le n$,
such that $\mu\not\subseteq\la$, \item $\I_\mu(\mu;\th;h)\ne0$.
\endroster

{\rm(ii)} These polynomials are unique up to scalar factors.
\endproclaim

Here the notation $\mu\not\subseteq\la$ means that the diagram
of $\mu$ is not contained in the diagram of $\la$.

\demo{Proof} Any polynomial in $x_1,\dots,x_n$ which is
$\Wa$--invariant in variables $x_i-\th i+ h$, $i=1,\dots,n$, can
be viewed simply as a symmetric polynomial in new variables
$(x_i-\th i+ h)^2$. Then the claims of the proposition become a
particular case of the results of \cite{Ok2}, corresponding to
the case of the ``perfect grid of class II''. \qed
\enddemo

Note that the existence of the polynomials
$\I_\mu(x_1,\dots,x_n;\th;h)$ can be derived from the results of
Okounkov \cite{Ok1} (see also Rains' paper \cite{R} which
contains a different approach to the results of \cite{Ok1}).
Namely,
$$
\I_\mu(x_1,\dots,x_n;\th;h)=\lim_{q\to1}
\frac{P^*_\mu(q^{x_1},\dots,q^{x_n};q,q^\th,q^{h-n\th})}
{(q-1)^{2|\mu|}}\,, \tag2.1
$$
where the polynomial in the numerator is the $BC_n$ type shifted
(or interpolation) Macdonald polynomial defined in \cite{Ok1,
Definition 1.3}.

The polynomials $\I_\mu(x_1,\dots,x_n;\th;h)$ are normalized by
setting
$$
\qquad\quad\I_\mu(\mu;\th;h)=\prod_{(i,j)\in\mu}
\big(1+\mu_i-j+\th(\mu'_j-i)\big)
\big(2h-1+\mu_i+j-\th(\mu'_j+i)\big)\,,\tag2.2
$$
where the product is over all squares $(i,j)$ in the diagram of
$\mu$. A motivation for such a normalization is given in
\cite{Ok2, Proposition 2.9}. Note that this normalization is
well adapted to the combinatorial formula \tht{2.4}. We shall
call the polynomials $\I_\mu(x_1,\dots,x_n;\th;h)$ the {\it
interpolation $BC_n$ polynomials\/}.

The next result provides a {\it combinatorial formula\/} for the
interpolation $BC_n$ polynomials. It is convenient to state it
in terms of {\it reverse\/} tableaux (cf. \cite{OO2}). A reverse
tableau $T$ of shape $\mu$ with entries in $\{1,\dots,n\}$ is
defined as a function $T(i,j)$ assigning to each box
$(i,j)\in\mu$ a number $T(i,j)\in\{1,\dots,n\}$, such that the
numbers decrease strictly down each column and weakly from left
to right along each row. The only difference with the
conventional definition of a semistandard tableau consists in
inverting the natural order in the set of indices $1,\dots,n$.

First, recall the combinatorial formula for the Jack polynomials
with parameter $\th$:
$$
P(x_1,\dots,x_n;\th)=\sum_T\psi_T(\th)
\prod_{(i,j)\in\mu}x_{T(i,j)}\,,\tag2.3
$$
summed over all reverse tableaux $T$ of shape $\mu$, as defined
above, where $\psi_T(\th)$ is a certain weight factor, which is
a rational function in $\th$ (see formulas (7.13'), (10.10),
(10.11), and (10.12) in chapter VI of Macdonald's book
\cite{Ma3}, and also \cite{OO2}; note that our parameter $\th$
is inverse to the parameter $\a$ used in \cite{Ma3}). We do not
need the explicit expression for $\psi_T(\th)$.

\proclaim{Proposition 2.2} We have
$$
\multline \I_\mu(x_1,\dots,x_n;\th;h)=\sum_T\psi_T(\th)
\prod_{(i,j)\in\mu} \big[ \big(x_{T(i,j)}+h-\th T(i,j)\big)^2\\
-\big((j-1)-\th (i-1)+h-\th T(i,j)\big)^2\big] \,,
\endmultline \tag2.4
$$
summed over all reverse tableaux $T$ of shape $\mu$, with
entries in $\{1,\dots,n\}$, where $\psi_T(\th)$ is the same
weight factor as in \tht{2.3}.
\endproclaim

\demo{Proof} Using \tht{2.1}, this can be obtained as a
degeneration of the combinatorial formula for the $BC_n$
interpolation Macdonald polynomials, established in Okounkov
\cite{Ok1, Theorem 5.2}. See also \cite{Ok2} and Rains \cite{R}.
\qed
\enddemo

Note that the expression for $I_\mu(\mu;\th;h)$ given in
\tht{2.2} can be obtained from the combinatorial formula
\tht{2.4}. Indeed, if $(x_1,\dots,x_n)=(\mu_1,\dots,\mu_n)$ then
all products in the right--hand side, except a single one
(corresponding to a special choice of $T$) vanish, and the only
nonvanishing term gives the expression \tht{2.2}. This can be
shown using the same argument as in the second proof of Theorem
11.1 from Okounkov--Olshanski \cite{OO1}.

We list a number of corollaries of these two propositions.

{}From the combinatorial formula \tht{2.4} it follows
immediately that the polynomials $\I_\mu(\,\cdot\,;\th;h)$,
which are initially defined as elements of $\Lt_n(h)$, actually
belong to $\Lt_n[h]$.

Next, the characterization of these polynomials given in
Proposition 2.1 shows that they are stable:
$$
\I_\mu(x_1,\dots,x_n,0;\th;h)= \I_\mu(x_1,\dots,x_n;\th;h),
$$
and, therefore, the sequence
$\{\I_\mu(x_1,\dots,x_n;\th;h)\}_{n\ge\ell(\mu)}$ correctly
defines an element of the algebra $\Lt[h]$. This element will be
denoted as $\I_\mu(x_1,x_2,\dots;\th;h)$ or
$\I_\mu(\,\cdot\,;\th;h)$ or else simply $I_\mu$.

Let us extend the filtration from $\Lt$ to $\Lt[h]$ by setting
$\deg h=1$. It is clear that the associated graded algebra
$\gr\Lt[h]$ is naturally isomorphic to the algebra $\La[h]$,
where $\La$ denotes the algebra of symmetric functions in
infinitely many variables.

Given an element $f\in\Lt[h]$, let us denote by $[f]$ its
highest degree term, which is a homogeneous element of the
algebra $\La[h]$. We remark that the element $\I_\mu\in\Lt[h]$
has degree $2|\mu|$ and it follows from the comparison of
\tht{2.3} and \tht{2.4} that
$$
[\I_\mu(x_1,x_2,\dots;\th;h)]
=P_\mu\big(x_1(x_1+2h),\,x_2(x_2+2h),\dots ;\th\big), \tag2.5
$$
where $P_\mu(\,\cdot\,;\th)\in\La$ is the Jack symmetric
function. Note that for any homogeneous symmetric function
$g(x_1,x_2,\dots)$ of degree $m$, the expression
$$
g(x_1(x_1+2h),\,x_2(x_2+2h),\dots)
$$
is a well defined homogeneous element of the algebra $\La[h]$ of
degree $2m$.

It is also clear that
$$
\I_\mu(\mu;\th;h)=H(\mu;\th)(2h)^{|\mu|}+\dots \tag2.6
$$
where the dots stand for the lower degree terms in $h$ and
$H(\mu;\th)$ is the following hook--length product
$$
H(\mu;\th):=\prod_{(i,j)\in\mu}((\mu_i-j)-\th(\mu'_j-i)+1) \,.
$$

\subhead 2B. Binomial formula
\endsubhead

Here by the binomial formula we mean an expansion of the
normalized Jacobi polynomials $\Phi_\la(z_1,\dots,z_n;\th,a,b)$
(see \tht{1.6}) about the point $(1,\dots,1)$.

\proclaim{Proposition 2.3} Let $\ell(\la)\le n$. We have
$$
\multline \Phi_\la(z_1,\dots,z_n;\th,a,b)\\
=\sum_{\mu} \frac{\I_\mu(\la;\th;\s+\th n) \,
P_\mu(z_1+z^{-1}_1-2,\dots, z_n+z^{-1}_n-2;\th)}
{C(n,\mu;\th;a,b)} \,. \endmultline\tag2.7
$$
where  $\I_\mu(\la;\th;\s+\th n)$ is the result of specializing
$h=\s+\th n$ in $\I_\mu(\la;\th;h)$, $\s=(a+b+1)/2$, and
$$
C(n,\mu;\th;a,b)= \I_\mu(\mu;\th;\s+\th n)\,
\J_\mu(\underbrace{1\dots,1}_n;\th,a,b).\tag2.8
$$
\endproclaim

By virtue of condition 3 in Proposition 2.1, the summation in
\tht{2.7} actually goes on the finite set of $\mu$'s such that
$\mu\subseteq\la$. An explicit expression for $C(n,\mu;\th;a,b)$
is given in Remark 2.5 below.

\demo{Proof} The expansion \tht{2.7} is a limit case of the
binomial formula for Koornwinder polynomials obtained in
\cite{Ok1, Theorem 7.1}. Namely, the Koornwinder polynomials are
orthogonal on $\T^n$ with weight
$$
\prod_{i<j\le n} \frac {(z^{\pm1}_i z^{\pm1}_j)_\infty}
{(tz^{\pm1}_i z^{\pm1}_j)_\infty} \prod_{i=1}^n
\frac{(z^{\pm1}_i, -z^{\pm1}_i,q^{1/2} z^{\pm1}_i, -q^{1/2}
z^{\pm1}_i)_\infty} {(a_1 z^{\pm1}_i, -a_2 z^{\pm1}_i,q^{1/2}
a_3 z^{\pm1}_i, -q^{1/2} a_4 z^{\pm1}_i)_\infty}\,,
$$
where $q,t,a_1,\dots,a_4$ are the 6 parameters and,
by definition,
$$
(u_1,u_2,\dots)_\infty=\prod_k \prod_{i=0}^\infty (1-q^i u_k)\,.
$$
\vskip  -3 \jot
\noindent If one sets
\vskip  -4 \jot
$$
t=q^\th\,, \quad a_1=q^{a+\frac12}\,, \quad a_2=q^{b+\frac12}\,,
\quad a_3=a_4=0\,,
$$
and lets $q\to 1$ then these 6-parametric polynomials become
$\J_\la(z;\th,a,b)$. Using this, \tht{2.1}, and the following
immediate corollary of Theorem 5.2 in \cite{Ok1}
$$
\lim_{q\to 1} P^*_\mu(z;q,q^\th,q^{(a+b+1)/2}) =
P_\mu(z_1+z^{-1}_1-2,\dots, z_n+z^{-1}_n-2;\th)
$$
one obtains \tht{2.7} from Theorem 7.1 in \cite{Ok1}. Here
$P^*_\mu$ stands for the 3-parametric $BC_n$-type interpolation
Macdonald polynomial defined in \cite{Ok1, Definition 1.3}; it
has already appeared in \tht{2.1}.

Alternatively, one can prove \tht{2.7} without going into
$q$-analogs by just repeating the proof of main theorem of
\cite{OO2} and using the fact that the algebra of commuting
differential operators whose eigenfunctions are $\J_\la$ is
isomorphic under the Harish--Chandra homomorphism to the algebra
of polynomials in $\la_1,\dots,\la_n$ that are $\Wa$-invariant
in variables
$$
\la+\rho=(\la_1 + (n-1)\th + \s,\dots,\la_{n-1}+\th+\s,\la_n+\s)
\,.
$$
\qed
\enddemo

Note that the expansion of Jacobi polynomials in Jack
polynomials was obtained much earlier by James and Constantine
\cite{JK} for the case when $\th=\frac12$ or 1. However, in
their formula the coefficients of the expansion are written in
quite a different form, as certain combinatorial sums over
standard tableaux of shape $\la/\mu$ (this is {\it not\/}
equivalent to combinatorial formula \tht{2.4}). Their result was
extended to arbitrary values of $\th$ by Macdonald (unpublished
work \cite{Ma1, \S9}; see also \cite{BeO, \S5}) and Lassalle
\cite{L}. A discussion of the role of binomial formulas for the
characters of the orthogonal and symplectic groups can be found
in \cite{OO3}.

\subhead 2C. Asymptotics of denominators in binomial formula
\endsubhead

We fix an arbitrary partition $\mu$ and let $n$ go to infinity.
As in Theorem 1.4, we assume that the parameters $a,b$ may
depend on $n$. We write them as $a_n,b_n$ and assume that the
limits \tht{1.7} exist.

\proclaim{Proposition 2.4} The denominator \tht{2.8} in
\tht{2.7} has the following asymptotics
$$
C(n,\mu;\th;a_n,b_n)\sim\frac{H(\mu;\th)}{H'(\mu;\th)}\,
4^{|\mu|}\,\th^{|\mu|}\,(\th+\bar a)^{|\mu|}\cdot
n^{2|\mu|}\,,\tag2.9
$$
where
$$
H(\mu;\th)=\prod_{(i,j)\in\mu}((\mu_i-j)-\th(\mu'_j-i)+1),
\qquad
H'(\mu;\th)=\prod_{(i,j)\in\mu}((\mu_i-j)-\th(\mu'_j-i)+\th)\,,
$$
and $\bar a=\lim a_n/n$ as in \tht{1.7}.
\endproclaim

\demo{Proof} Recall that $C(n,\mu;\th;a_n,b_n)$ is the product
of two terms, $\I_\mu(\mu;\th;\s_n+\th n)$ and
$\J_\mu(1,\dots,1;\th,a_n,b_n)$, where $\s_n=(a_n+b_n+1)/2$. We
claim that, as $n\to\infty$, the following two asymptotic
relations hold
$$
\gather
 \I_\mu(\mu;\th;\s_n+\th n)\,\sim\, H(\mu;\th)\,2^{|\mu|}\,(\th+\bar\s)^{|\mu|}
n^{|\mu|}\,, \tag2.10 \\
\J_\mu(\underbrace{1,\dots,1}_n;\th,a_n,b_n)\, \sim \,
\frac1{H'(\mu;\th)}\,2^{|\mu|}\,\left(\frac{\th+\bar
a}{\th+\bar\s}\right)^{|\mu|}\,\th^{|\mu|}\, n^{|\mu|}\,,
\tag2.11
\endgather
$$
where
$$
\bar\s=\lim_{n\to\infty}\frac{\s_n}n=\frac{\bar a+\bar b}2\,.
$$

Clearly, \tht{2.10} and \tht{2.11} imply \tht{2.9}.

The first relation immediately follows from \tht{2.6}, let us
check the second relation.

The following is the general formula, due to Opdam, for the
value of a multivariate Jacobi polynomial, indexed by a weight
$\mu$, at the unit element, see \cite{HS}, Part I, Theorem
3.6.6,
$$
\prod_{\a>0} \frac{\G\left((\mu+\rho,\a^\vee)+k_\a+\frac12
k_{\a/2}\right)} {\G\left((\mu+\rho,\a^\vee)+\frac12
k_{\a/2}\right)} \frac{\G\left((\rho,\a^\vee)+\frac12
k_{\a/2}\right)} {\G\left((\rho,\a^\vee)+k_\a+\frac12
k_{\a/2}\right)}\,,
$$
where $\a^\vee$ stands for the root dual to $\a$, and
$k_{\a/2}=0$ if the root $\a/2$ does not exist.

In our case, the polynomial in question is just
$\J_\mu(\,\cdot\,;\th,a_n,b_n)$, and the unit element is
identified with the point $(1,\dots,1)$. Next, we have
$$
\rho=((n-1)\th + \s_n,\dots,\th+\s_n,\s_n)
$$
and there are 4 types of the positive roots $\a$
$$
\e_i-\e_j\,, \quad \e_i+\e_j \quad (1\le i<j\le n), \qquad
\e_i\,, \quad 2\e_i \quad (1\le i\le n)
$$
with formal multiplicities
$$
k_{\e_i\pm\e_j}=2\th, \quad k_{\e_i}=a_n-b_n\,, \quad
k_{2\e_i}=b_n+\tfrac12\,.
$$
As the scalar product used in Opdam's formula we may take the
natural scalar product in $\R^n$. Then the dual roots $\a^\vee$
are as follows
$$
(\e_i\pm\e_j)^\vee=\e_i\pm\e_j\,, \quad \e_i^\vee=2\e_i\,, \quad
(2\e_i)^\vee=\e_i\,.
$$

We split the product over $\a>0$ into 4 products according to
these 4 types of positive roots
$$
\J_\mu(\underbrace{1,\dots,1}_{\text{$n$ times}};\th,a_n,b_n)
=\ppm \ppp \pp \pdp \,, \tag2.12
$$
where
$$
\alignat2 &\ppm &&= \prod_{1\le i<j\le n}
\frac{\G(\mu_i-\mu_j+\th(j-i+1))} {\G(\mu_i-\mu_j+\th(j-i))}
\frac{\G(\th(j-i))}
{\G\th(j-i+1))} \,, \\
&\ppp &&=\prod_{1\le i<j\le n}
\frac{\G(\mu_i+\mu_j+\th(2n-i-j+1)+2\s_n)}
{\G(\mu_i+\mu_j+\th(2n-i-j)+2\s_n)} \frac{\G(\th(2n-i-j)+2\s_n)}
{\G(\th(2n-i-j+1)+2\s_n)}\,,\\
&\pp &&= \prod_{1\le i \le n} \frac{\G(2\mu_i+2\th(n-i)+2a_n+1)}
{\G(2\mu_i+2\th(n-i)+2\s_n)} \frac{\G(2\th(n-i)+2\s_n)}
{\G(2\th(n-i)+2a_n+1)}\,,\\
&\pdp &&= \prod_{1\le i \le n} \frac{\G(\mu_i+\th(n-i)+2\s_n)}
{\G(\mu_i+\th(n-i)+a_n+\tfrac12)}
\frac{\G(\th(n-i)+a_n+\tfrac12)} {\G(\th(n-i)+2\s_n)}\,.
\endalignat
$$

It is clear that $\ppm$ is just the formula for the value of the
Jack polynomial (which is essentially the Jacobi polynomial for
the $A$ series). Therefore (see Stanley \cite{St} or Macdonald
\cite{Ma3}, VI.10.12)
$$
\ppm\sim H'(\mu;\th)^{-1} \th^{|\mu|}\, n^{|\mu|}\,.
$$
Next, we have
$$
\gather \pp\sim\left( \frac{\th n + a_n}{\th n
+\s_n}\right)^{2|\mu|}\sim\,
\left(\frac{\th + \bar a}{\th +\bar\s}\right)^{2|\mu|},\\
\pdp\sim\left( \frac{\th n + 2\s_n}{\th n +
a_n}\right)^{|\mu|}\, \sim\left( \frac{\th + 2\bar\s}{\th + \bar
a}\right)^{|\mu|}.
\endgather
$$
Finally, consider the product $\ppp$. All factors with
$\ell(\mu)<i<j$ are trivial (equal to 1). Moreover, the
contribution of every fixed pair $i<j$ is asymptotically
trivial, so that the whole contribution of the factors with
$1\le i<j\le\ell(\mu)$ is trivial. Therefore,
$$
\ppp\sim \prod\Sb 1\le i\le\ell(\mu)\\ \ell(\mu)<j\le n\endSb
\frac{\G(\mu_i+\th(2n-i-j+1)+2\s_n)}
{\G(\mu_i+\th(2n-i-j)+2\s_n)} \frac{\G(\th(2n-i-j)+2\s_n)}
{\G(\th(2n-i-j+1)+2\s_n)}\,.
$$
For every $i$ the product over $j$ telescopes to
$$
\gather \frac{\G(\mu_i+\th(2n-i-\ell(\mu))+2\s_n)}
{\G(\mu_i+\th(n-i)+2\s_n)} \frac{\G(\th(n-i)+2\s_n)}
{\G(\mu_i+\th(2n-i-\ell(\mu))+2\s_n)} \\
\sim \left( \frac{2\th n +2\s_n}{\th n + 2\s_n}\right)^{\mu_i}
\, \sim \left( \frac{2\th  +2\bar\s}{\th +
2\bar\s}\right)^{\mu_i}\,,
\endgather
$$
so that
$$
\ppp\sim 2^{|\mu|}\left( \frac{\th  +\bar\s}{\th +
2\bar\s}\right)^{|\mu|}.
$$

It follows that
$$
\ppp \pp \pdp \sim 2^{|\mu|}\,\left(\frac{\th+\bar
a}{\th+\bar\s}\right)^{|\mu|}.
$$
This gives formula \tht{2.11} and concludes the proof. \qed
\enddemo

\remark{Remark 2.5} Using \tht{2.12} one can derive the
following explicit expression for the quantity \tht{2.8}
$$
C(n,\mu;\th;a,b)=4^{|\mu|}\,\frac{H(\mu;\th)}{H'(\mu;\th)}\,
\prod_{i=1}^n\frac{\Gamma(\mu_i+(n-i+1)\th)\Gamma(\mu_i+(n-i)\th+a+1)}
{\Gamma((n-i+1)\th)\Gamma((n-i)\th+a+1)}\,.
$$
In particular, $C(n,\mu;\th;a,b)$ actually does not depend on
$b$. The asymptotic relation \tht{2.9} is readily obtained from
this expression.
\endremark

\head
3. Sufficient conditions of regularity
\endhead

Recall that we have fixed some $\th>0$. Define elements
$g_1,g_2, \dots$ in the algebra $\La$ of symmetric functions by
means of a generating function
$$
1+\sum_{k=1}^\infty g_k(x_1,x_2,\dots)t^k
=\prod_{j=1}^\infty(1-x_j t)^{-\th}\,,
$$
where $t$ is a formal variable. These elements are algebraically
independent generators of $\La$. Given some VK parameters
$\a,\be,\ga$, define an algebra homomorphism
$$
\eabg: \La \longrightarrow \C
$$
by
$$
1+\sum_{k=1}^\infty g_k t^k @>\quad\eabg\quad>> e^{\ga\th
t}\prod_{i=1}^\infty \frac{1+\be_i\th t}{(1-\a_i t)^\th}
$$
(cf. the definition of the extended symmetric functions, see
\cite{OO4, Section 2.7} or \cite{KOO}). In a less formal way,
this homomorphism can be written as
$$
\prod_{j=1}^\infty(1-x_j t)^{-\th} @>\quad\eabg\quad>> e^{\ga\th
t}\prod_{i=1}^\infty \frac{1+\be_i\th t}{(1-\a_i t)^\th}\,.
$$
We will also use such a notation in the sequel.

\proclaim{Proposition 3.1} Let $\{\la(n)\}$ be a VK sequence of
partitions with parameters $\a,\be,\ga$, see Definition 1.3, and
let $\{h_n\}$ be a sequence of complex numbers such that the
limit
$$
\bar h=\lim_{n\to\infty}\frac{h_n}n
$$
exists. Then for any element $f\in\Lt[h]$ we have
$$
\lim_{n\to\infty}\frac{f(\la(n);h_n)}{n^{\deg f}}
=\eabg\left(\big[f\big]\Big|_{h=\bar h}\right)\,,
$$
where $\big[f\big]\in\La[h]$ is the highest degree term of $f$.
\endproclaim

\demo{Proof} Assume first that $f$ does not depend on $h$, that
is, $f\in\Lt$. Then $\big[ f\big]\in\La$, and the claim is that
$$
\lim_{n\to\infty}\frac{f(\la(n))}{n^{\deg f}}
=\eabg\left(\big[f\big]\right)\,,
$$
which follows from Theorem 7.1 in \cite{KOO} (see also Theorem
3.1 in \cite{OO4}). Returning to the general case, we write
$f\in\Lt[h]$ as a polynomial in $h$,
$$
f=f_0+f_1h+f_2h^2+\dots+f_mh^m, \qquad f_i\in\Lt,
$$
and apply the above formula to each $f_0,\dots,f_m$. \qed
\enddemo

Given $\tau\in\C$, define a homomorphism
$$
\pta: \La \to \La
$$
by the formula
$$
(\pta f)(x_1,x_2,\dots)\mapsto
f(x_1(x_1+\tau),x_2(x_2+\tau),\dots)\,, \qquad f\in\La\,.
\tag3.1
$$
Then we have the following

\proclaim{Proposition 3.2}
$$
\eabg\,\pta \left(\prod_{j=1}^\infty\frac1{(1-x_j t)^{\th}}
\right)= e^{\ga \th\tau t}\prod_{i=1}^\infty
\frac{1+\th\be_i(\tau-\th\be_i)t}{(1-\a_i(\tau+\a_i)t)^\th}\,.
$$
\endproclaim
\demo{Proof}
By the definition of $\pta$ and $\eabg$ we have
$$
\align \eabg\pta \left(\prod_{i=1}^\infty\frac1{(1-x_j t)^{\th}}
\right)&=\eabg\left(\prod_{j=1}^\infty\frac1{(1-x_j(x_j+\tau)
t)^{\th}}
\right)\\
&= \eabg\left(\prod_{j=1}^\infty\frac1{((1-x_j t_1)(1-x_j
t_2))^{\th}}
\right)\,,\\
\intertext{where $t_1+t_2=\tau t$ and $t_1t_2=-t$\,,} &=e^{\ga
\th(t_1+t_2)} \prod_{i=1}^\infty \frac{1+\be_i\th t_1}{(1-\a_i
t_1)^\th}
\frac{1+\be_i\th t_2}{(1-\a_i t_2)^\th}\\
&=e^{\ga \th \tau t}\prod_{i=1}^\infty
\frac{1+\th\be_i(\tau-\th\be_i)t}{(1-\a_i(\tau+\a_i)t)^\th}\,,
\endalign
$$
which is the result stated. \qed
\enddemo

Now we can prove the following

\proclaim{Theorem 3.3}  Assume that the condition \tht{1.7} is
fulfilled. Let $\{\la(n)\}$ be a Vershik--Kerov sequence with
parameters $\a,\be,\ga$, see Definition 1.3. For any fixed
$k=1,2,\dots$, the functions
$$
\Phi_{\la(n)} (z_1,\dots,z_k,\underbrace{1,\dots,1
}_{\text{$n-k$ times}};\th,a_n,b_n)
 \tag3.2
$$
converge infinitesimally, in the sense of Definition 1.2, to the
function
$$
\prod_{l=1}^k \phi_{\a,\be,\ga,\bar a,\bar
b}\left(\frac{z_l+z^{-1}_l}2\right)\,, \tag3.3
$$
where $\bar a=\lim a_n/n$, $\bar b=\lim b_n/n$ as in \tht{1.7},
and $\phi_{\a,\be,\ga,\bar a,\bar b}$ is the following function
of a single variable $x\in[-1,1]$
$$
\phi_{\a,\be,\ga,\bar a,\bar b}(x)=e^{\ga(x-1)}
\prod_{i=1}^{\infty} \dfrac{
1+\dfrac{\be_i}2\left(\dfrac{2\th+\bar a+\bar
b-\th\be_i}{\th+\bar a}\right)(x-1)}
{\left(1-\dfrac{\a_i}{2\th}\left(\dfrac{2\th+\bar a+\bar
b+\a_i}{\th+\bar a}\right) (x-1)\right)^\th}\,.
$$
\endproclaim

\demo{Proof} Apply to the functions \tht{3.2} the binomial
formula \tht{2.7}. In that formula $\mu$ ranges over partitions
with $\ell(\,\cdot\,)\le n$. However, in our case, the last
$n-k$ arguments in $P_\mu(\,\cdot\,;\th)$ equal 0. By the
stability of Jack polynomials, this implies that the summation
actually goes over partitions $\mu$ with $\ell(\mu)\le k$. Thus,
we obtain
$$
\multline \Phi_{\la(n)} (z_1,\dots,z_k,\underbrace{1,\dots,1
}_{\text{$n-k$ times}};\th,a_n,b_n)=\\
=\sum_{\mu, \ell(\mu)\le k\,} \frac{\I_\mu(\la(n);\th;\s_n+\th
n) \, 2^{|\mu|}} {C(n,\mu;\th,a_n,b_n)} \,
P_\mu(x_1-1,\dots,x_k-1;\th)
\endmultline \tag3.4
$$
where
$$
x_l=\frac{z_l+z^{-1}_l}2\,, \quad l=1,\dots,k\,.
$$
Since the polynomials $P_\mu(\,\cdot\,;\th)$ form a homogeneous
basis in the space of symmetric polynomials,  the infinitesimal
convergence of the left--hand side of \tht{3.4} in the sense of
Definition 1.2, as $n\to\infty$, is equivalent to the
coefficient--wise convergence of the expansion in the
right--hand side. Thus, we have to examine the asymptotics of
the quantities
$$
\frac{\I_\mu(\la(n);\th;\s_n+\th n) \, 2^{|\mu|}}
{C(n,\mu;\th;a_n,b_n)}, \qquad \text{$\mu$ fixed, $n\to\infty$.}
$$

By Proposition 2.4,
$$
C(n,\mu;\th;a_n,b_n)\sim\frac{H(\mu;\th)}{H'(\mu;\th)}\,
4^{|\mu|}\,\th^{|\mu|}\,(\th+\bar a)^{|\mu|}\cdot n^{2|\mu|}\,,
$$
Next, by Proposition 3.1  we have
$$
\lim_{n\to\infty} \frac{\I_\mu(\la(n);\s_n+\th n)} {n^{2|\mu|}}
= \ep_{\a,\be,\ga}
\left(\big[\I_\mu(\,\cdot\,;h)\big]\Big|_{h=\th+\bar\s}\right)
\,,
$$
and by \tht{2.5} and \tht{3.1}
$$
\big[\I_\mu(\,\cdot\,;h)\big]\Big|_{h=\th+\bar\s}=\pi_{2(\th+\bar\s)}
P_\mu(\,\cdot\,;\th),
$$
which implies
$$
\lim_{n\to\infty} \frac{\I_\mu(\la(n);\s_n+\th n)} {n^{2|\mu|}}
= \ep_{\a,\be,\ga}\,\pi_{2\th+2\bar\s}\, P_\mu(\,\cdot\,;\th).
$$
Therefore, the expansion \tht{3.4} converges coefficient--wise
to
$$
\e_{\a,\be,\ga} \pi_{2(\th+\bar\s)} \left( \sum_{\mu}
\frac{H'(\mu)}{H(\mu)}\, P_\mu(y_1,y_2,\dots;\th)\,
P_\mu\left(\frac{x_1-1}{2\th(\th+\bar a)},
\dots,\frac{x_k-1}{2\th(\th+\bar a)};\th\right)\right)\,,
\tag3.5
$$
where the homomorphisms $\e_{\a,\be,\ga}$ and
$\pi_{2(\th+\bar\s)}$ act on the variables $y_1,y_2,\dots$.  The
sum in the round brackets equals
$$
\multline \sum_{\mu} Q_\mu(y_1,y_2,\dots;\th)\,
P_\mu\left(\frac{x_1-1}{2\th(\th+\bar a)},
\dots,\frac{x_k-1}{2\th(\th+\bar a)};\th\right)\\
=\prod_{l=1}^k \prod_{j=1}^\infty \left(1-\frac{y_j
(x_l-1)}{\th(\th+\bar a)}\right)^{-\th}.
\endmultline
$$
Here we have used the Cauchy identity for Jack polynomials (see
Macdonald \cite{Ma3, ch. VI, \tht{4.13} and \S10}). Therefore
\tht{3.5} equals
$$
\prod_{l=1}^k \e_{\a,\be,\ga} \pi_{2\th+2\bar\s}
\left(\prod_{j=1}^\infty \left(1-\frac{y_j (x_l-1)}{\th(\th+\bar
a)}\right)^{-\th}\right)\,.
$$
Using Proposition 3.2 and the relation $2\bar\s=\bar a+\bar b$
one transforms this into the desired result. \qed
\enddemo

\proclaim{Corollary 3.4} In the hypotheses of Theorem 3.3,
assume additionally that $a_n\ge b_n\ge-\frac12$. Then for any
$k=1,2,\dots$, the functions \tht{3.2} converge to the function
\tht{3.3} uniformly on the torus $\T^k$.
\endproclaim

\demo{Proof} The assumption $a_n\ge b_n\ge-\frac12$ makes it
possible to apply Proposition 1.1 which implies that the
functions \tht{2.6} are positive definite functions on the torus
$\T^k$. Since the function \tht{3.3} is real--analytic, our
claim follows from Theorem 3.3 by virtue of a well--known
general fact (see, e.g., Lemma 4.2 from \cite{OO4}). \qed
\enddemo

\head 4. Necessary conditions of regularity
\endhead
The argument of this section is similar to that of \cite{OO4,
\S5}. We begin with three technical lemmas which then are used
to prove Proposition 4.4.

\proclaim{Lemma 4.1} Assume $h \ge \th n - 1/2$. Then we have
$$
\I_{(2)}(\la_1,\dots,\la_n;\th;h) \le
\left(\I_{(1)}(\la_1,\dots,\la_n;\th;h)\right)^2
$$
for any partition $\la$, $n=1,2,\dots$, and $\th\ge 0$.
\endproclaim

Observe that if $h=\th n+ (a+b+1)/2$, where $a,b>-1$, then the
hypothesis of the lemma is satisfied.

\demo{Proof} Set $A=\I_{(1)}(\la;\th;h)$ and
$B=\I_{(2)}(\la;\th;h)$. From the combinatorial formula
\tht{2.4} we obtain
$$
\align
A&=\sum_i \left((\l_i+h-\th i)^2 - (h-\th i)^2\right)\,,\\
B&=\sum_i \left((\l_i+h-\th i)^2 - (h-\th i)^2\right)
\left((\la_i+h-\th i)^2 - (h-\th i+1)^2\right)\\
&\phantom{=}+\frac{2\th}{1+\th} \sum_{i<j} \left((\l_j+h-\th
j)^2 - (h-\th j)^2\right) \left((\la_i+h-\th i)^2 - (h-\th
i+1)^2\right)\,.
\endalign
$$
It is elementary to check that for any $l\in\Z_{\ge 0}$ and any
$i\le n$ one has
$$
(l+h-\th i)^2 - (h-\th i)^2 \ge 0
$$
provided $h \ge \th n - 1/2$; in particular,
$$
(1+h-\th i)^2 - (h-\th i)^2 \ge 0 \,.
$$
Therefore
$$
\align B&\le\sum_i \left((\la_i+h-\th i)^2 - (h-\th i)^2\right)
\left((\la_i+h-\th i)^2 - (h-\th i)^2\right)\\
&\phantom{\le}+\frac{2\th}{1+\th} \sum_{i<j} \left((\la_j+h-\th
j)^2 - (h-\th j)^2\right) \left((\la_i+h-\th i)^2 - (h-\th
i)^2\right)
\endalign
$$
(here we use the inequality with general $l$ to conclude that
the first factors are nonnegative, while the inequality with
$l=1$ is used to remove 1 in the second factors).  Since
$\frac{2\th}{1+\th} < 2$ for any $\th > -1$, the lemma follows.
\qed
\enddemo

\proclaim{Lemma 4.2} Suppose that
$$
\I_{(1)}(\la(n);\th;h_n)=O(n^2)\,, \quad n\to\infty \,,
$$
and
$$
h_n\sim h_0 n\,, \quad h_0 > \th/2, \quad n\to\infty\,.
$$
Then
$$
|\la(n)|=O(n)\,, \quad n\to\infty\,.
$$
\endproclaim

Observe that if $h_n=\s_n+\th n$, where $\s_n\sim\bar\s n$,
$\bar\s\ge 0$, then the hypothesis of the lemma is satisfied.

\demo{Proof} We have
$$
\align \I_{(1)}(\la;\th;h_n)&=
\sum_{i=1}^n \la_i(\la_i+2h_n-2\th i)\\
&=\sum_{i=1}^n \la_i^2+\th\sum_{1\le i<j\le n}
(\la_i-\la_j)+(2h_n-\th(n+1)) \sum_{i=1}^n \la_i \,.
\endalign
$$
Since all summands are positive for large $n$ we conclude
that
$$
(2h_0-\th)\, n \sum_{i=1}^n \la_i(n) = O(n^2)\,, \quad
n\to\infty\,,
$$
which implies $|\la(n)|=O(n)$. \qed
\enddemo

\proclaim{Lemma 4.3} Let $\bar a\ge0$ and $\bar b\ge0$ be fixed.
Assume $(\dd\a, \dd\be, \dd\ga)$ and $(\ddd\a, \ddd\be,
\ddd\ga)$ are two systems of VK parameters such that the
corresponding functions \tht{1.11} coincide,
$$
\gather e^{\dd{\ga} (x-1)} \prod_{i=1}^{\infty} \dfrac{
1+\dfrac{\dd{\be}_i}2\left(\dfrac{2\th+\bar a+\bar
b-\th\dd{\be}_i}{\th+\bar a}\right) (x-1)}
{\left(1-\dfrac{\dd{\a}_i}{2\th}\left(\dfrac{2\th+\bar a+\bar
b+\dd{\a}_i}{\th+\bar a}\right) (x-1)\right)^\th}\\
= e^{\ddd{\ga} (x-1)} \prod_{i=1}^{\infty} \dfrac{
1+\dfrac{\ddd{\be}_i}2\left(\dfrac{2\th+\bar a+\bar
b-\th\ddd{\be}_i}{\th+\bar a}\right) (x-1)}
{\left(1-\dfrac{\ddd{\a}_i}{2\th}\left(\dfrac{2\th+\bar a+\bar
b+\ddd{\a}_i}{\th+\bar a}\right) (x-1)\right)^\th}
\endgather
$$
for any $x\in[-1,1]$.

Then $(\dd\a, \dd\be, \dd\ga)=(\ddd\a, \ddd\be, \ddd\ga)$.
\endproclaim

\demo{Proof} A similar claim was established in step 3 of the
proof of Theorem 5.1 from \cite{OO4}. We shall use a similar
argument.

Write the above identity of functions as
$$
e^{\dd{\ga} y} \prod_{i=1}^{\infty} \dfrac{ 1+\dd{B_i}y}
{\left(1-\dd{A_i}y\right)^\th}=e^{\ddd{\ga} y}
\prod_{i=1}^{\infty} \dfrac{ 1+\ddd{B_i}y}
{\left(1-\ddd{A_i}y\right)^\th}\,, \qquad y\in[-2,0],  \tag4.1
$$
and observe that
$$
\gather \dd{A_1}\ge\dd{A_2}\ge\dots\ge0, \qquad
\dd{B_1}\ge\dd{B_2}\ge\dots\ge0,
\qquad\sum_{i=1}^\infty(\dd{A}_i+\dd{B_i})<\infty \\
\ddd{A_1}\ge\ddd{A_2}\ge\dots\ge0, \qquad
\ddd{B_1}\ge\ddd{B_2}\ge\dots\ge0,
\qquad\sum_{i=1}^\infty(\ddd{A}_i+\ddd{B_i})<\infty
\endgather
$$
Since the correspondences
$$
\dd{\a}_i\leftrightarrow\dd{A}_i, \qquad
\dd{\be}_i\leftrightarrow\dd{B}_i, \qquad
\ddd{\a}_i\leftrightarrow\ddd{A}_i, \qquad
\ddd{\be}_i\leftrightarrow\ddd{B}_i
$$
are one--to--one, it suffices to prove that \tht{4.1} implies
$$
\dd{A}_i=\ddd{A}_i, \qquad \dd{B}_i=\ddd{B}_i, \qquad
\dd{\ga}=\ddd{\ga}.
$$

We may extend \tht{4.1} to an identity between two holomorphic
functions in $y$. These functions are well defined at least in
the left half--plane $\Re y<0$. Actually, the left--hand side is
holomorphic in the half--plane $\Re y<(\dd{A}_1)^{-1}$ and has a
singularity at $y=(\dd{A}_1)^{-1}$. Likewise, the right--hand
side is holomorphic in the half--plane $\Re y<(\ddd{A}_1)^{-1}$
and has a singularity at $y=(\ddd{A}_1)^{-1}$. This implies
$\dd{A}_1=\ddd{A}_1$. Thus, both sides of \tht{4.1} have common
factors which can be cancelled. Iterating this procedure we
prove that $\dd{A}_2=\ddd{A}_2$, etc. Then we come to an
identity of entire functions,
$$
e^{\dd{\ga} y} \prod_{i=1}^{\infty} (1+\dd{B_i}y) =e^{\ddd{\ga}
y} \prod_{i=1}^{\infty} (1+\ddd{B_i}y).
$$
Examining the zeros of both sides we see that
$\dd{B}_i=\ddd{B}_i$ for all $i=1,2,\dots$, and finally we
conclude that $\dd{\ga}=\ddd{\ga}$. \qed
\enddemo

\proclaim{Theorem 4.4} Let $a_n\ge b_n\ge-\frac12$ and let
$\{\la(n)\}$ be a sequence of partitions with $\ell(\la(n))\le
n$. Assume that either $\{\la(n)\}$ is minimally regular, or the
functions
$$
\Phi_{\la(n)}(z,\underbrace{1,\dots,1 }_{\text{$n-1$
times}};\th,a_n,b_n)\,, \quad |z|=1\,, \tag4.2
$$
converge infinitesimally about $z=1\in\T^1$, see Definition 1.2.

Then $\{\la(n)\}$ is a Vershik-Kerov sequence.
\endproclaim

\demo{Proof} The proof given below is completely parallel to the
proof of Theorem 5.1 in \cite{OO4}.

{\it Step 1\/}. Consider the following function on the unit
circle $\T$
$$
\phi_n(z)=\Phi_{\la(n)}(z,\,\underbrace{1,\dots,1}_{n-1}\,
;\th,a_n,b_n), \qquad z\in\T.
$$
By the binomial formula \tht{2.7} and Proposition 2.4 we have
$$
\phi_n(z)=1+A_{1,n}(z+z^{-1}-2)+A_{2,n}(z+z^{-1}-2)^2+\dots,
$$
where
$$
A_{1,n}\sim\const_1\,\frac{\I_{(1)}(\la(n);\th n+\s_n)}{n^2}\,,
\qquad A_{2,n}\sim\const_2\,\frac{\I_{(2)}(\la(n);\th
n+\s_n)}{n^4}\,, \tag4.3
$$
with some constants not depending on $n$.

We claim that
$$
\I_{(1)}(\la(n);\th n+\s_n)=O(n^2). \tag4.4
$$

Indeed, in case when the functions \tht{4.2} converge
infinitesimally, this bound is immediate. Let us prove it when
$\{\la(n)\}$ is minimally regular. If $z=e^{i\varphi}$ then
$$
z+z^{-1}-2=2(\cos\varphi-1)=-\varphi^2+\tfrac1{12}\varphi^4+O(\varphi^6)
$$
so that
$$
\gather \phi_n(e^{i\varphi})
=1+A_{1,n}(-\varphi^2+\tfrac1{12}\varphi^4)+A_{2,n}\varphi^4+O(\varphi^6)\\
=1-A_{1,n}\varphi^2+(A_{2,n}+\tfrac1{12}A_{1,n})\varphi^4+O(\varphi^6).
\endgather
$$

On the other hand, we know that $\phi_n(z)$ is a normalized
positive definite function on $\T$ (see the proof of Corollary
3.4), hence it is the Fourier transform of a probability measure
$M_n$ on the lattice $\Z$. The assumption of minimal regularity
means that the measures $M_n$ weakly converge to a probability
measure on $\Z$. In such a situation, Lemma 5.2 from \cite{OO4}
says that if the second moments of the measures $M_n$ are not
uniformly bounded then, passing to a suitable subsequence where
the second moments tend to infinity, we obtain that the fourth
moments grow faster than the squared second moments. Up to
constant factors, the second and fourth moments are the
coefficients in $\varphi^2$ and $\varphi^4$, respectively. This
means that if the number sequence $\{A_{1,n}\}$ is unbounded
then, for a suitable subsequence of indices $n$, $A_{2,n}$ grows
faster than $(A_{1,n})^2$. But, by virtue of \tht{4.3}, this
contradicts Lemma 4.1. We conclude that $\{A_{1,n}\}$ is
bounded, which is equivalent to \tht{4.4}.

{\it Step 2\/}. We claim that
$$
|\la(n)|=O(n)\,, \quad n\to\infty\,. \tag4.5
$$
Indeed, this follows from the result of step 1 and Lemma 4.2.

{\it Step 3\/}. Using the bound \tht{4.5} and Cantor's diagonal
process, we see that any subsequence of $\{\la(n)\}$ contains a
VK subsequence. It remains to prove that any two VK subsequences
of $\{\la(n)\}$ have the same VK parameters. By Theorem 3.3 and
the hypotheses of Theorem 4.4, both subsequences lead to one and
the same limit function. Then we apply Lemma 4.3 to conclude
that the VK parameters are the same. \qed
\enddemo

\demo{Proof of Theorem 1.4} Let, as usual, $\{\la(n)\}$ be a
sequence of partitions with $\ell(\la(n))\le n$, and let the
parameters $a_n,\,b_n$ satisfy the inequality $a_n\ge
b_n\ge-\frac12$ and condition \tht{1.7}.

By Theorem 4.4, each of the 3 regularity properties of
Definition 1.2 implies the VK property of Definition 1.3.
Conversely, assume $\{\la(n)\}$ is a VK sequence. Then, by
virtue of Theorem 3.3 and Corollary 3.4, this implies all 3
regularity properties of Definition 1.2. This proves claim (i)
of Theorem 1.4. Corollary 3.4 also proves claim (ii) of Theorem
1.4. \qed
\enddemo

\head 5. The convex set $\Ut$ \endhead

In this section we fix parameters $a,b$ such that $a\ge
b\ge-\frac12$.

For $n=1,2,\dots$ let $\U_n$ be the set of functions on the
torus $\T^n$ of the form
$$
\varphi(z_1,\dots,z_n)=\sum_{\la:\,\ell(\la)\le n}
c_\la\Phi_\la(z_1,\dots,z_n;\th,a,b),
$$
where
$$
c_\la\ge0, \qquad \sum_{\la:\,\ell(\la)\le n} c_\la=1.
$$
Note that the coefficients $c_\la$ are uniquely determined by
the function $\varphi$, because the Jacobi polynomials form an
orthogonal basis in a suitable $L^2$ space. Recall that each
$\Phi_\la(\,\cdot\,;\th,a,b)$ is a positive definite function on
the torus (see the proof of Corollary 3.4), normalized at the
unit element $(1,\dots,1)\in\T^n$. It follows that
$|\Phi_\la(z_1,\dots,z_n;\th,a,b)|\le1$, which implies that the
series converges uniformly and defines a continuous function on
$\T^n$. Thus, $\U_n$ is a subset of the set of continuous,
positive definite, normalized functions on $\T^n$. It is clear
that $\U_n$ is a convex set. As an abstract convex set, it is
isomorphic to a simplex with infinitely many vertices.

\proclaim{Proposition 5.1} Let $\la$ be an arbitrary partition
with with $\ell(\la)\le n$. Expand the function
$\J_\la(z_1,\dots,z_{n-1},1;\th,a,b)$ in Jacobi polynomials in
$n-1$ variables $z_1,\dots,z_{n-1}$, with the same parameters
$\th,a,b$. Then all coefficients in this expansion are
nonnegative.
\endproclaim

\demo{Proof} This fact can be derived by a degeneration from the
branching rule for Koornwinder polynomials, established by Rains
\cite{R, (5.76)}. \qed
\enddemo

By Proposition 5.1, the specialization $z_n=1$ sends
$\Phi_\la(z_1,\dots,z_n;\th,a,b)$ to a function from $\U_{n-1}$.
Hence this specialization map determines an affine map
$\U_n\to\U_{n-1}$. Using these maps for $n=2,3,\dots$ we set
$$
\Ut=\varprojlim \U_n\,, \qquad n\to\infty.
$$
As we shall see (Corollary 5.3), this projective limit space
does not depend on $a,b$.

An equivalent definition is as follows. Consider the set
$\T^\infty_0=\varinjlim\T^n$ whose elements are infinite vectors
$(z_1,z_2,\dots)\in\T\times\T\times\dots$ with finitely many
coordinates $z_i$ distinct from 1. Then elements of $\Ut$  can
be described as functions $\varphi(z_1,z_2,\dots)$ on $\T^\infty_0$
such that for any $n=1,2,\dots$, the function
$$
\varphi_n(z_1,\dots,z_n)=\varphi(z_1,\dots,z_n,1,1,\dots)
$$
on $\T^n$ belongs to $\U_n$.

It is clear that $\Ut$ is a convex set. Let $\Ex\Ut$ denote the
set of its extreme points.

\proclaim{Theorem 5.2} There is a one--to--one correspondence
between elements of the set $\Ex\Ut$ and collections
$(\a,\be,\ga)$ of VK parameters,
$$
\a=(\a_1\ge\a_2\ge\dots\ge0), \quad
\be=(\be_1\ge\be_2\ge\dots\ge0), \quad \ga\ge0, \quad
\sum(\a_i+\be_i)<\infty.
$$

Given $(\a,\be,\ga)$, the corresponding element of $\Ex\Ut$,
viewed as a function on $\T^\infty_0$, has the form
$$
\Phi_{\a,\be,\ga}(z_1,z_2,\dots)
=\prod_{j=1}^\infty\phi_{\a,\be,\ga}\left(\frac{z_j+z_j^{-1}}2\right)
$$
where $(z_1,z_2,\dots)\in\T^\infty_0$, and the function
$\phi_{\a,\be,\ga}$ is as in Theorem 1.4,
$$
\phi_{\a,\be,\ga}\left(\frac{z+z^{-1}}2\right)=e^{\frac\gamma2(z+z^{-1}-2)}
\prod_{i=1}^\infty \frac{\left(1+\frac12\beta_i(z-1)\right)
\left(1+\frac12\beta_i(z^{-1}-1)\right)}
{\left(\left(1-\frac12\alpha_i(z-1)/\th\right)
\left(1-\frac12\alpha_i(z^{-1}-1)/\th\right)\right)^\theta}
$$
or equivalently
$$
\phi_{\a,\be,\ga}(x) =e^{\gamma (x-1)}\,\prod_{i=1}^\infty
\dfrac{\left(1+\beta_i\left(1-\dfrac{\be_i}2\right)\right)(x-1)}
{\left\{\left(1-\dfrac{\a_i}{\th}\left(1+\dfrac{\a_i}{2\th}\right)\right)
(x-1)\right\}^\theta}\,, \qquad x\in[-1,1].
$$

In particular, $\Ex\Ut$ does not depend on the parameters $a,b$.
\endproclaim

\demo{Proof} This can be proved exactly as Theorem 1.4 in
\cite{OO4}, see \S6 in \cite{OO4}. By virtue of a general result
(see Theorem 6.1 in \cite{OO4}), each function $\Phi\in\Ex\Ut$
can be approximated (uniformly on any finite--dimensional
sub--torus $\T^k$) by a sequence
$\{\Phi_{\la(n)}(\,\cdot\,;\th,a,b)\}$. By Theorem 1.4, the
limit functions
$$
\Phi=\lim_{n\to\infty}\Phi_{\la(n)}(\,\cdot\,;\th,a,b)
$$
are precisely the functions of the form $\Phi_{\a,\be,\ga}$.
This shows that $\Ex\Ut\subset\{\Phi_{\a,\be,\ga}\}$.

The inverse inclusion $\{\Phi_{\a,\be,\ga}\}\subset\Ex\Ut$ is
obtained by a simple argument using de Finetti's theorem.
Indeed, the fact that each function $\Phi_{\a,\be,\ga}$ is
positive definite on $\T^\infty_0$ and has the multiplicative
form $\phi(z_1)\phi(z_2)\dots$ implies that $\Phi_{\a,\be,\ga}$
is an extreme point in a larger convex set, namely the set of
characteristic functions of symmetric probability measures on
$\Z^\infty=\Z\times\Z\times\dots$.  \qed
\enddemo

\proclaim{Proposition 5.3} Any element of the convex set $\Ut$
is represented by a probability measure on the set $\Ex\Ut$ of
extreme points, and this representation is unique.
\endproclaim

\demo{Proof} This claim does not follow directly from Choquet's
theorem because the set $\Ut$ is not compact. However, it can be
checked, for instance, by the method of \cite{O4, \S9}.\qed
\enddemo

\proclaim{Corollary 5.4} The set $\Ut$ does not depend on the
parameters $a,b$.
\endproclaim

\demo{Proof} Indeed, according to Theorem 5.2, $\Ex\Ut$ does not
depend on $a,b$. Then the corollary follows from Proposition
5.3. \qed
\enddemo

\head 6. Spherical functions on infinite--dimensional symmetric
spaces \endhead

By an infinite--dimensional symmetric space we mean a
homogeneous space $G/K$ where the $G$ and $K$ are inductive
limits of groups,
$$
G=\varinjlim G(n), \quad K=\varinjlim K(n), \qquad n\to\infty,
$$
such that for any $n$, $G(n)/K(n)$ is a Riemannian symmetric
space of compact type belonging to one of the classical series.
We assume that $\rank (G(n)/K(n))=n$.

There are 10 such spaces $G/K$ corresponding to 10 classical
series of symmetric spaces, see Olshanski \cite{O1}, \cite{O2}:

\centerline{Table I}

\medskip

1. $G(n)=U(n)$, $K(n)=O(n)$.

2. $G(n)=U(n)\times U(n)$, $K(n)=U(n)$.

3. $G(n)=U(2n)$, $K(n)=Sp(n)$.

4. $G(n)=O(2n)$, $K(n)=O(n)\times O(n)$.

5. $G(n)=Sp(n)$, $K(n)=U(n)$.

6. $G(n)=U(2n)$, $K(n)=U(n)\times U(n)$.

7. $G(n)=O(\wt n)\times O(\wt n)$, $K(n)=O(\wt n)$.

8. $G(n)=Sp(n)\times Sp(n)$, $K(n)=Sp(n)$.

9. $G(n)=Sp(2n)$, $K=Sp(n)\times Sp(n)$.

10. $G(n)=O(2\wt n)$, $K(n)=U(\wt n)$.

\medskip

\remark{Comments} a) The embeddings $G(n)\to G(n+1)$ and
$K(n)\to K(n+1)$ which are implicit in the definition of the
groups $G$ and $K$ are natural ones. The embeddings $K(n)\to
G(n)$ are also quite evident. In particular, for series
2,\,7,\,8, these are the diagonal embeddings.

b) For series 7 and 10, we wrote $\wt n$ instead of $n$ because
in these cases the rank equals $[\wt n/2]$. Here one may choose
one of the two possible variants: $\wt n=2n$ or $\wt n=2n+1$,
and the inductive limit space $G/K$ does not depend of the
choice, up to isomorphism.

c) For the Grassmann spaces (series 4,\,6,\,9), we could equally
well use two distinct indices $n_1$, $n_2$. That is, we could
deal with the spaces $O(n_1+n_2)/O(n_1)\times O(n_2)$,
$U(n_1+n_2)/U(n_1)\times U(n_2)$, and $Sp(n_1+n_2)/Sp(n_1)\times
Sp(n_2)$. Again, such a generalization does affect the limit
space $G/K$, provided that both indices go to infinity.

d) In the case of series 2,\,7,\,8, the symmetric space
$G(n)/K(n)$ is one of the classical groups $U(n)$, $O(n)$,
$Sp(n)$, and the corresponding infinite--dimensional space $G/K$
coincides with one of groups
$$
U(\infty)=\varinjlim U(n), \qquad O(\infty)=\varinjlim O(n),
\qquad Sp(\infty)=\varinjlim Sp(n).
$$

\endremark

As shown in \cite{O1}, \cite{O2}, for the 10 pairs $(G,K)$
listed in Table I there is a rich theory of unitary
representations. The present paper concerns a part of this
theory related to spherical representations.

Assume $(G,K)$ is one of the 10 pairs from Table I. Let $T$ be a
unitary representation of $G$ in a Hilbert space $H$, and
$\xi\in H$ be a distinguished unit $K$--invariant vector. We say
that $(T,\xi)$ is a {\it spherical representation\/} of the pair
$(G,K)$ if $\xi$ is a cyclic vector. That is, if the linear span
of the vectors $T(g)\xi$, where $g$ ranges over $G$, is dense in
$H$. Two spherical representations, $(T_1,\xi_1)$ and
$(T_2,\xi_2)$, are said to be {\it equivalent\/} if there is an
isometry $H_1\to H_2$ of the corresponding Hilbert spaces taking
$\xi_1$ to $\xi_2$ and commuting with the action of $G$.

Attached to any spherical representation $(T,\xi)$ is its {\it
spherical function\/}
$$
F(g)=(T(g)\xi,\xi), \qquad g\in G.
$$
This is a positive definite function on the group $G$,
two--sided invariant with respect to the subgroup $K$, and
taking value 1 at the unity $e\in G$. Denote by $\Up(G,K)$ the
set of all functions with these three properties.  Then the
correspondence $T\to F$ defined above determines a {\it
bijection\/} between equivalence classes of spherical
representations and functions from the set $\Up(G,K)$.

Of special interest are {\it irreducible\/} spherical
representations. That is, those $(T,\xi)$ for which $T$ is an
irreducible unitary representation of $G$. If $T$ is an
irreducible unitary representation of $G$ then a $K$--invariant
vector $\xi$ (provided it exists) is {\it unique\/}, within a
scalar factor which does not affect the spherical function. The
spherical functions of irreducible spherical representations are
precisely the {\it extreme points\/} of $\Up(G,K)$ (it is worth
noting that $\Up(G,K)$ is a convex set).

Thus, classifying the irreducible spherical representations of
$(G,K)$ is equivalent to describing the extreme points of the
convex set $\Up(G,K)$. We aim to explain how this problem is
related to that discussed in \S6.

Let $\Up(G(n),K(n))$ denote the set of functions on $G(n)$ that
are positive definite, two--sided $K(n)$--invariant and take
value 1 at the unity. This is a convex set, isomorphic to an
infinite--dimensional simplex. Note that the vertices of this
simplex, which are the extreme points of $\Up(G(n),K(n))$, are
the {\it indecomposable spherical functions\/} of $(G(n),K(n))$,
that is, matrix coefficients of the form
$$
(T^{(n)}(g)\xi^{(n)},\xi^{(n)}), \qquad g\in G(n),
$$
where $T^{(n)}$ is an arbitrary irreducible finite--dimensional
unitary representation of $G(n)$ possessing a unit
$K(n)$--invariant vector $\xi^{(n)}$ (such a vector is unique,
within a scalar factor of absolute value 1).

The natural embedding of pairs
$$
(G(n-1),K(n-1))\hookrightarrow (G(n),K(n))
$$
induces an affine projection
$$
\Up(G(n),K(n))\;\to\;\Up(G(n-1),K(n-1))
$$
and we have
$$
\Up(G,K)=\varprojlim \Up(G(n),K(n)), \qquad n\to\infty.
$$

Below we focus on the series 4--10, the case of series 1--3
being the subject of our previous paper \cite{OO4}.

So, let $(G(n),K(n))$ belong to one of the seven series 4--10
from Table I. Let $\Cal R_n$ denote the restricted root system
of the symmetric space $(G(n),K(n))$. Then $\Cal R_n$ coincides
with one of the classical root system $B_n$, $C_n$, $D_n$ or
$BC_n$ equipped with appropriate root multiplicities.

In all cases, it is convenient to regard $\Cal R_n$ as a
subsystem of the $BC_n$ root system
$$
R_n=\left\{\pm\e_i\pm\e_j: 1\le i<j\le
n\right\}\cup\left\{\pm\e_i,\pm2\e_i: 1\le i\le n\right\}.
$$
To each of the series 4--10 one can attach a particular triple
of parameters $\th,a,b$ in such a way that the formal
multiplicities $k_\a$ defined in \tht{1.5} coincide with the
true root multiplicities in $\Cal R_n$. In particular, $k_\a=0$
means that the root $\a\in R_n$ is not contained in $\Cal R_n$.

Specifically, we have:
\medskip

\centerline{Table II}

\medskip
4. $G(n)=O(2n)$, $K(n)=O(n)\times O(n)$:\quad $\th=\tfrac12$,
$a=b=-\tfrac12$, $\Cal R_n=D_n$.

5. $G(n)=Sp(n)$, $K(n)=U(n)$:\quad $\th=\tfrac12$, $a=b=0$,
$\Cal R_n=C_n$.

6. $G(n)=U(2n)$, $K(n)=U(n)\times U(n)$:\quad $\th=1$, $a=b=0$,
$\Cal R_n=C_n$.

7. $G(n)=O(\wt n)\times O(\wt n)$, $K(n)=O(\wt n)$:\quad
$\th=1$, $a=-\frac12$ or $\frac12$, $b=-\tfrac12$, $\Cal
R_n=D_n$ or $B_n$.

8. $G(n)=Sp(n)\times Sp(n)$, $K(n)=Sp(n)$:\quad $\th=1$,
$a=b=\tfrac12$, $\Cal R_n=C_n$.

9. $G(n)=Sp(2n)$, $K=Sp(n)\times Sp(n)$:\quad $\th=2$, $a=b=1$,
$\Cal R_n=C_n$.

10. $G(n)=O(2\wt n)$, $K(n)=U(\wt n)$:\quad $\th=2$, $a=0$ or 2,
$b=0$, $\Cal R_n=C_n$ or $BC_n$.

\medskip

Here, in the case of series 7 and 10, the first option for $a$
and $\Cal R_n$ is chosen  if $\wt n=2n$, and the second option
is chosen if $\wt n=2n+1$.

\proclaim{Proposition 6.1} Let $(G(n),K(n))$ belong to one of
the series 4--10 and let $\th,a,b$ be the corresponding
parameters as listed above. Assume additionally that the
restricted root system $\Cal R_n$ is not $D_n$. Then there are
natural bijections
$$
\U_n\longleftrightarrow \Up(G(n),K(n)), \qquad n=1,2,\dots
\tag6.1
$$
which are isomorphisms of convex sets and commute with the
projections
$$
\U_n\to\U_{n-1}\,, \qquad \Up(G(n),K(n))\to \Up(G(n-1),K(n-1)).
\tag6.2
$$
\endproclaim

\demo{Idea of proof} Let $\Cal W_n$ denote the restricted Weyl
group of the symmetric space $G(n)/K(n)$. For any classical
series, $\Cal W_n$ may be identified with the $BC_n$ Weyl group
$S(n)\ltimes\Z^2_n$ (realized as the group $W_*$ acting on the
torus $\T^n$ as explained in \S1) or with its subgroup of index
2 (the latter possibility holds exactly when $\Cal R_n=D_n$).

It is well known that the indecomposable spherical functions of
$(G(n),K(n))$ can be interpreted as normalized $\Cal
W_n$--invariant orthogonal polynomials on the torus $\T^n$ with
the weight \tht{1.1}, where the parameters $\th,a,b$ are those
attached to the corresponding series. When $\Cal R_n\ne D_n$,
these are precisely the normalized polynomials
$\Phi_\la(z_1,\dots,z_n;\th,a,b)$.

This yields the required bijection \tht{6.1}. The fact that
these bijections are compatible with the projections \tht{6.2}
is readily verified. \qed
\enddemo

\proclaim{Corollary 6.2} Under the hypotheses of Proposition 6.1
we have: the convex set $\Up(G,K)$ is isomorphic to the convex
set $\Ut$ described in \S5. In particular, it depends only on
$\th$ but not on $a,b$.
\endproclaim

\demo{Proof} Indeed, this follows at once from Proposition 6.1
and the results of $\S5$. \qed
\enddemo

\remark{Remark 6.3} When the restricted root system $\Cal R_n$
is of type $D_n$, certain polynomials
$\Phi_\la(z_1,\dots,z_n;\th,a,b)$ turn out to be half--sums of
two distinct $\Cal W_n$--invariant orthogonal polynomials. As a
consequence, the set $\Up(G(n),K(n))$ turns out to be somewhat
larger than the set $\U_n$. Nevertheless, the claim of Corollary
6.2 holds in this case as well. Indeed, as is seen from the list
above, the equality $\Cal R_n=D_n$ occurs in two cases: for
series 4 and for series 7 with $\wt n$ even. In the latter case
we may choose $\wt n$ odd without changing the limit space
$G/K$. In the former case the same effect is achieved if we take
$G(n)=O(2n+m)$, $K(n)=O(n+m)\times O(n)$, where $m$ is an
arbitrary fixed positive integer.
\endremark

\head 7. The $BC_n$ polynomials with $\th=1$ \endhead

When $\th=1$, the multivariate Jacobi polynomials $\J_\la$ and
the interpolation polynomials $\I_\mu$ admit explicit
determinantal expressions. This makes it possible to establish
the basic facts about these polynomials independently of the
general theory, in a rather elementary way.

Throughout the present section we assume $\th=1$ and fix
arbitrary parameters $a>-1,b>-1$. We give explicit formulas for
both kinds of polynomials and sketch elementary proofs of the
binomial formula \tht{2.7} and of the branching rule for Jacobi
polynomials, which in turns implies Proposition 5.1.

Let $\p_l(x;a,b)$ (where $l=0,1,2,\dots$) denote the classical
Jacobi polynomials in a single variable $x$, orthogonal on the
segment $-1\le x\le 1$ with the weight function
$(1-x)^a(1+x)^b$. We use the same normalization as in Erdelyi et
al. \cite{Er}. An explicit expression for $\p_l(x;a,b)$ in terms
of the Gauss hypergeometric function ${}_2F_1$ is
$$
\p_l(x;a,b)=\frac{\G(l+a+1))}{\G(l+1)\G(a+1)}\,
{}_2F_1(-l,l+a+b+1;a+1;\tfrac{1-x}2). \tag7.1
$$
In particular, the value at $x=1$ is given by
$$
\p_l(1;a,b)=\frac{\G(l+a+1))}{\G(l+1)\G(a+1)}
$$
and the leading coefficient in $\p_l(x;a,b)$ is
$$
\varkappa(l;a,b)=2^{-l}\,\frac{\G(2l+a+b+1)}{\G(l+a+b+1)\G(l+1)}\,.
$$

More generally, for any $n=1,2,\dots$ and any partition $\la$
with $\ell(\la)\le n$ we set
$$
\PP_\la(x_1,\dots,x_n;a,b) =\frac{\det\limits_{1\le i,j\le
n}\left[\p_{\la_i+n-i}(x_j\,;a,b)\right]}{V(x)}\,,  \tag7.2
$$
where
$$
V(x)=\prod_{1\le i<j\le n}(x_i-x_j).
$$

\proclaim{Proposition 7.1} The $BC_n$ orthogonal polynomials
with $\th=1$ and arbitrary parameters $a,b>-1$ are expressed
through the polynomials \tht{7.2} as follows
$$
\J_\la(z_1,\dots,z_n;1,a,b) =\const\,
\PP_\la\left(\dfrac{z_1+z_1^{-1}}2\,,\dots,\dfrac{z_n+z_n^{-1}}2\,;a,b\right)
$$
where
$$
\const=\frac{2^{|\la|}}{\prod_{i=1}^n\varkappa(\la_i+n-i;a,b)}
$$
\endproclaim

This fact is undoubtedly well known. E.g., in an equivalent
form, it was pointed out in Lassalle \cite{L}. For reader's
convenience we present a proof.

\demo{Proof} It is readily verified that the polynomials
\tht{7.2} are pairwise orthogonal on the $n$--dimensional cube
$[-1,1]^n$ with respect to the measure
$$
V^2(x) \prod_{1\le i \le n} (1-x_i)^a (1+x_i)^b \, dx_1\dots,
dx_n
$$
which implies that the polynomials
$$
(z_1,\dots,z_n)\mapsto
\PP_\la\left(\dfrac{z_1+z_1^{-1}}2\,,\dots,\dfrac{z_n+z_n^{-1}}2\,;a,b\right)
$$
are pairwise orthogonal on the torus $\T^n$ with the weight
\tht{1.1} specialized at $\th=1$.

Next, we have
$$
\PP_\la(x_1,\dots,x_n;a,b) = \const' \, x_1^{\la_1}\dots
x_n^{\la_n} +\dots
$$
where dots mean lower terms in lexicographic order and
$$
\const'=\prod_{i=1}^n\varkappa(\la_i+n-i;a,b),
$$
which implies the triangularity condition \tht{1.2} for the
polynomials on the torus, defined by the right--hand side of
\tht{7.2}.

The stronger triangularity condition \tht{1.4} can also be
readily verified. Indeed, comparing the determinantal expression
\tht{7.2} with the determinantal formula for the Schur
polynomials
$$
s_\mu(x_1,\dots,x_n)=\frac{\det\limits_{1\le i,j\le
n}\left[x_j^{\mu_i+n-i}\right]}{V(x)}\,,
$$
we see that
$$
\PP_\la(x_1,\dots,x_n;a,b)=\sum_{\mu\subseteq\la}a_{\la\mu}s_\mu(x_1,\dots,x_n)
$$
where $\mu\subseteq\la$ means that the diagram of $\mu$ is
contained in that of $\la$, and $a_{\la\mu}$ are certain
coefficients. Next, recall the well--known triangularity
property of the Schur polynomials:
$$
s_\mu(x_1,\dots,x_n)=\sum_{\nu\le\mu}K_{\mu\nu}m_\nu(x_1,\dots,x_n),
$$
where $m_\nu$ is the conventional monomial symmetric function,
$K_{\mu\nu}$ are certain coefficients (the Kostka numbers) and
$\nu\le\mu$ is the dominance order on partitions ($\mu-\nu$ can
be written as a linear combination of the vectors $\e_i-\e_j$,
$i<j$).

{}From the last two formulas we obtain the triangularity
condition of the form \tht{1.4}:
$$
\gather
\PP_\la\left(\dfrac{z_1+z_1^{-1}}2\,,\dots,\dfrac{z_n+z_n^{-1}}2\,;a,b\right)
=\sum_{\nu\ll\la}b_{\la\nu}
m_\nu\left(\dfrac{z_1+z_1^{-1}}2\,,\dots,\dfrac{z_n+z_n^{-1}}2\right)\\
=\sum_{\mu\ll\la}c_{\la\mu} \wt m_\nu(z_1\,,\dots,z_n)
\endgather
$$
with certain coefficients $b_{\la\nu}$ and $c_{\la\mu}$.

Thus, the polynomials in $(z_1,\dots,z_n)$ defined by the
right--hand of \tht{7.2} possess  the characteristic properties
of the $BC_n$ orthogonal polynomials on the torus $\T^n$ with
parameter $\th=1$ and hence coincide with the polynomials
$\J_\la(z_1,\dots,z_n;1,a,b)$.

\qed
\enddemo

By virtue of \tht{7.2}, Proposition 7.1 provides an explicit
determinantal expression for the polynomials
$\J_\la(z_1,\dots,z_n;1,a,b)$. Now we shall give an explicit
expression of the $BC_n$ interpolation polynomials with $\th=1$.

We need a notation. Given an infinite sequence of parameters
$A=(A_1,A_2,\dots)$, define ``generalized powers'' of a variable
$y$ by
$$
(y\mid A)^m=(y-A_1)\dots(y-A_m), \quad m=1,2,\dots; \qquad
(y\mid A)^0=1.
$$
Next, we define (generalized) {\it factorial Schur
polynomials\/} in $n$ variables by
$$
s_\mu(y_1,\dots,y_n\mid A)=\frac{\det\limits_{1\le i,j\le
n}\left[\,(y_i\mid A)^{\mu_j+n-j}\,\right]}{V(y_1,\dots,y_n)}\,,
\tag7.3
$$
where $\mu$ is an arbitrary partition with $\ell(\mu)\le n$.
Note that $s_\mu(y_1,\dots,y_n\mid A)$ is an inhomogeneous
symmetric polynomial whose top degree homogeneous component is
the conventional Schur polynomial $s_\mu(y_1,\dots,y_n)$. The
polynomials \tht{7.3} share many properties of the conventional
Schur polynomials. A number of formulas for the polynomials
\tht{7.3} can be found in Macdonald \cite{Ma2} and \cite{Ma3,
Example I.3.20}, Molev \cite{Mo}, and also, for the special case
$A=(0,1,2,\dots)$, in Okounkov--Olshanski \cite{OO1}. Note that
our notation for the parameters $A_1,A_2,\dots$ differs from
that of Macdonald by sign.

\proclaim{Proposition 7.2} The $BC_n$ interpolation polynomials
with $\th=1$ are expressed through the factorial Schur
polynomials \tht{7.3} as follows
$$
\I_\mu(x_1,\dots,x_n;1;h)=s_\mu((x_1+h-1)^2,\dots,(x_n+h-n)^2\mid
A), \tag7.1
$$
where
$$
A=((h-n)^2,(h-n+1)^2,(h-n+2)^2,\dots).
$$
\endproclaim

\demo{Proof} We have to check that in the special case $\th=1$,
all claims of Proposition 2.1 hold and the corresponding
polynomials coincide with those given by formula \tht{7.3}. This
was shown, in an elementary way, in Okounkov--Olshanski
\cite{OO3, Theorem 2.5}. \qed
\enddemo

About this result, see also \cite{Ok2, \S3.3}.

\proclaim{Proposition 7.3} The polynomials
$I_\mu(x_1,\dots,x_n;1;h)$ defined by \tht{7.1} satisfy the
combinatorial formula of Proposition 2.2. That is,
$$
\multline \I_\mu(x_1,\dots,x_n;1;h)=\sum_T
\prod_{(i,j)\in\mu} \big[ \big(x_{T(i,j)}+h- T(i,j)\big)^2\\
-\big(j-i+h- T(i,j)\big)^2\big] \,,
\endmultline
$$
summed over all reverse tableaux $T$ of shape $\mu$, with
entries in $\{1,\dots,n\}$.
\endproclaim

\demo{Proof} This is a special case of the combinatorial formula
for factorial Schur polynomials,
$$
s_\mu(y_1,\dots,y_n\mid A)=\sum_T \prod_{(i,j)\in\mu}
(y_{T(i,j)}-A_{j-i+n+1-T(i,j)}),
$$
summed over all reverse tableaux $T$ of shape $\mu$, with
entries in $\{1,\dots,n\}$. About the latter formula, see
Goulden--Greene \cite{GG}, Macdonald \cite{Ma2}, and also
Okounkov \cite{Ok2, \S3.3}. \qed
\enddemo

The binomial formula \tht{2.7} of Proposition 2.3 reduces to the
following claim

\proclaim{Proposition 7.4} Let $\ell(\la)\le n$. We have
$$
\frac{\PP_\la(x_1,\dots,x_n;a,b)}
{\PP_\la(\,\underbrace{1,\dots,1}_n\,;a,b)} =\sum_{\mu}
\frac{I_\mu(\la;1;\s+n) \, s_\mu(x_1-1,\dots, x_n-1)}
{c(n,\mu;a)} \,.
$$
where  $I_\mu(\la;1;\s+n)$ is the result of specializing
$h=\s+n$ in $I_\mu(\la;1;h)$, $\s=(a+b+1)/2$, and
$$
c(n,\mu;a)= 2^{|\mu|}\, \prod_{i=1}^n
\frac{\Gamma(\mu_i+n-i+1)\Gamma(\mu_i+n-i+a+1)}
{\Gamma(n-i+1)\Gamma(n-i+a+1)}
$$
\endproclaim

\demo{Proof} By virtue of the definition of the polynomials
$I_\mu(\,\cdot\,;1,h)$ this formula can be rewritten as
$$
\multline
\frac{\PP_\la(x_1,\dots,x_n;a,b)}{\PP_\la(\,\underbrace{1,\dots,1}_n\,;a,b)}\\
=\sum_{\mu} \frac{s_\mu(l_1^2,\dots,l^2_n\mid
\s^2,(\s+1)^2,(\s+2)^2,\dots) \, s_\mu(x_1-1,\dots, x_n-1)}
{c(n,\mu;a)} \,
\endmultline
$$
where
$$
l_i=\la_i+n-i+\s, \qquad i=1,\dots,n.
$$
The latter formula can be directly derived from the
determinantal formula \tht{7.2} and the expression \tht{7.1} for
the classical Jacobi polynomials. For some special values of the
parameters $a,b$ (which correspond to characters of classical
groups of the $B$, $C$, $D$ series) such a computation was done
in Okounkov--Olshanski \cite{OO3, Theorem 1.2}. For general
$a,b$ the argument is quite similar. \qed
\enddemo

The next result is a special case of Proposition 5.1.

\proclaim{Proposition 7.5} Let $n=2,3,\dots$ and $\la$ be a
partition with $\ell(\la)\le n$. In the expansion
$$
\PP_\la(x_1,\dots,x_{n-1},1;a,b)=\sum_{\nu:\;\ell(\nu)\le
n-1}(\dots)\,\PP_\nu(x_1,\dots,x_{n-1};a,b)
$$
all coefficients $(\dots)$ are nonnegative.
\endproclaim

\demo{Proof} a) First, let us describe the scheme of the proof.
We shall use certain renormalized polynomials
$$
\gathered R_\la(x_1,\dots,x_n)=\text{(a positive factor)}\,
\PP_\la(x_1,\dots,x_n;a,b)\\
\wt R_\mu(x_1,\dots,x_{n-1})=\text{(a positive factor)}\,
\PP_\la(x_1,\dots,x_{n-1};a+1,b).
\endgathered \tag7.4
$$
Below $\mu$ and $\nu$ denote partitions with $\ell(\,\cdot\,)\le
n-1$. We shall establish the following two--step branching rule,
which implies the proposition:
$$
\gather
R_\la(x_1,\dots,x_{n-1},1) =\sum_{\mu\prec\la}\wt R_\mu(x_1,\dots,x_{n-1}) \tag7.5\\
\wt R_\mu(x_1,\dots,x_{n-1})
=\sum_{\nu\prec\mu\cup0}A(\mu,\nu)\,R_\nu(x_1,\dots,x_{n-1}),
\qquad A(\mu,\nu)>0, \tag7.6
\endgather
$$
where $\mu\prec\la$ means
$$
\la_1\ge\mu_1\ge\la_2\ge\dots\ge\la_{n-1}\ge\mu_{n-1}\ge\la_n
$$
and $\nu\prec\mu\cup0$ means
$$
\mu_1\ge\nu_1\ge\mu_2\ge\dots\ge\mu_{n-1}\ge\nu_{n-1}\ge0.
$$
The coefficients $A(\mu,\nu)$ are specified below, see
\tht{7.13} and \tht{7.14}.

The proof presented below is quite elementary (we found it in
1998). It is worth noting that a much more general two--step
branching rule, in the context of Koornwinder polynomials, was
established by Rains, see \cite{R, (5.76)}.

b) Now we proceed to the proof of the proposition. For
$l=0,1,2,\dots$ set
$$
r_l(x)=\frac{\p_l(x;a,b)}{\p_l(1;a,b)}, \qquad \wt
r_l(x)=\frac{r_{l+1}(x)-r_l(x)}{x-1}\,. \tag7.7
$$
We have
$$
\gathered r_l(x)=\text{(a positive factor)}\,\p_l(x;a,b),\\
\wt r_l(x)=\text{(a positive factor)}\,\p_l(x;a+1,b).
\endgathered \tag7.8
$$
Indeed, the first relation in \tht{7.8} is evident, because
$\p_l(1;a,b)>0$. Let us check the second relation.

Since $r_{l+1}(1)=r_l(1)=1$, it is clear that $\wt r_l(x)$ is a
polynomial, and its degree is strictly equal to $l$. We have
$$
\gather
\int_{-1}^1 x^m \wt r_l(x)(1-x)^{a+1}(1+x)^b dx\\
=-\int_{-1}^1 x^m r_{l+1}(x)(1-x)^a(1+x)^b dx+\int_{-1}^1 x^m
r_l(x)(1-x)^a(1+x)^bdx,
\endgather
$$
and the latter two integrals vanish whenever $m<n$. This implies
that $\wt r_l(x)$ is proportional to $\p_l(x;a+1,b)$.

By the definition of $\wt r_l(x)$, its leading coefficient is
the same as that of $r_{l+1}(x)$, hence positive. This completes
the proof of \tht{7.8}.

c) Set
$$
\gather R_\la(x_1,\dots,x_n)=\frac{\det\limits_{1\le i,j\le
n}\left[r_{\la_i+n-i}(x_j)\right]}{V(x_1,\dots,x_n)} \tag7.9\\
\wt R_\mu(x_1,\dots,x_{n-1})=\frac{\det\limits_{1\le i,j\le
n-1}\left[\wt r_{\mu_i+n-1-i}(x_j)\right]}{V(x_1,\dots,x_{n-1})}
\tag7.10
\endgather
$$
Then, due to \tht{7.8}, we have \tht{7.4}.

d) Let us check the first step of our branching rule, \tht{7.5}.
By virtue of \tht{7.9}, the left--hand side of \tht{7.5} is
given by the ratio of a determinant of the order $n$ and a
Vandermonde. Examine the determinant in the numerator. For
$i=1,\dots,n-1$, let us subtract the $i+1$th row from the $i$th
row. Then we come to a determinant of order $n-1$. Its $(i,j)$th
entry is equal to
$$
r_{\la_i+n-i}(x_j)-r_{\la_{i+1}+n-1-i}(x_j),
$$
and we may divide it by $(x_j-1)$, because of the obvious
relation
$$
V(x_1,\dots,x_{n-1},1)=V(x_1,\dots,x_{n-1})\prod_{j=1}^{n-1}(x_j-1).
$$
Then the $(i,j)$ entry will take the form
$$
\frac{r_{l_1+1}(x_j)-r_{l_2}(x_j)}{x_j-1}= \sum_{m=l_2}^{l_1}\wt
r_m(x_j),
$$
where we abbreviated $l_1=\la_i+n-i-1$, $l_2=\la_{i+1}+n-1-i$.

Employing the latter expression and expanding the determinant
along the rows we get the desired result \tht{7.5}.

e) Let us check the second step of the branching rule,
\tht{7.6}.

Write the three--term recurrence relation for the polynomials
$r_m$:
$$
r_{m+1}(x)=(a_mx+b_m)r_m(x)-c_m r_{m-1}(x), \qquad m\ge1.
\tag7.11
$$
By the normalization,
$$
a_m+b_m-c_m=1.
$$
By making use of these relations and the definition of $\wt r_m$
we get
$$
\wt r_m=a_m r_m+c_m \wt r_{m-1}, \qquad m\ge1.
$$
Iterating this relation we further get
$$
\wt r_m=\sum_{l=k}^m B(m,l) r_l+ c_m\dots c_k \wt r_{k-1},
\qquad m\ge k\ge0, \tag7.12
$$
where
$$
B(m,l)=\left(\prod_{l<p\le m}c_p\right)a_l\,, \qquad m\ge l\ge0.
$$
When $k=0$ we agree that $\wt r_{-1}=0$ so that the last term in
\tht{7.12} disappears.

By virtue of \tht{7.9} and \tht{7.10}, the relation \tht{7.6} is
equivalent to
$$
\det\limits_{1\le i,j\le n-1}[\wt r_{\mu_i+n-1-i}(x_j)]
=\sum_{\nu\prec\mu\cup0}A(\mu,\nu)\,\det\limits_{1\le i,j\le
n-1} [r_{\nu_i+n-1-i}(x_j)]
$$
Examine the determinant in the left--hand side. Its entries in
the $i$th row are of the form $\wt r_{\mu_i+n-1-i}(x_j)$, where
$j=1,\dots, n-1$. We shall apply to them the relation
\tht{7.12}, taking $m=\mu_i+n-1-i$, $k=\mu_{i+1}+n-1-i$.

First, do this for the first row, $i=1$. We get a decomposition
of the form
$$
\wt
r_{\mu_1+n-2}(x_j)=\sum_{l=\mu_2+n-2}^{\mu_1+n-2}(...)r_l(x_j)+
(...)\wt r_{\mu_2+n-3}(x_j),
$$
where the coefficients marked as $(...)$ are expressed through
the $c$-- and $a$--coefficients. Remark that the last term
coincides, within a scalar factor, with that in the second row.
Consequently, when we expand the determinant along the first
row, it will play no role.

Now, we perform this expansion and then look at the second row
and repeat the same procedure, etc. For the $i$th row
($i=1,\dots,n-2$) we get the decomposition
$$
\wt r_{\mu_i+n-1-i}(x_j)
=\sum_{l=\mu_{i+1}+n-1-i}^{\mu_i+n-1-i}(...)r_l(x_j)+ (...)\wt
r_{\mu_{i+1}+n-2-i}(x_j),
$$
and when we come to the last row ($i=n-1$) then we choose $k=0$
so that the ``last term'' mentioned above will disappear at all.
Finally we get the desired expression \tht{7.6} with
coefficients $A(\mu,\nu)$ given by
$$
A(\mu,\nu)=\prod_{i=1}^{n-1}B(\mu_i+n-1-i,\, \nu_i+n-1-i)\,.
\tag7.13
$$

e) To check that these coefficients are strictly positive we use
a well--known general property of orthogonal polynomials: the
coefficients $a_m$, $c_m$ in the three--term relation \tht{7.11}
are strictly positive provided that the leading coefficients of
the polynomials are strictly positive (see \cite{Er}, section
10.3, formulas (7) and (8)). The latter property holds in our
case, so that we conclude that $B(m,l)>0$ for any $m\ge l\ge0$,
and finally $A(\mu,\nu)>0$.

f) The coefficients $B(m,l)$ entering the formula \tht{7.13} can
be explicitly computed:
$$
\multline B(m,l)\\=\frac{(2m+a+b)\Gamma(m+b+1)m!
(2l+a+b+1)\Gamma(l+a+b+1)\Gamma(l+a+1)}
{2\Gamma(m+a+b+2)\Gamma(m+a+2)\Gamma(n+b+1)l!}\,.
\endmultline\tag7.14
$$
When $l=0$, the product $(2l+a+b+1)\Gamma(l+a+b+1)$ must be
replaced by $\Gamma(l+a+b+2)=\Gamma(a+b+2)$.

Indeed, let $k_m$ be the leading coefficients in $\frak
p_m(x;a,b)$ and $h_m$ be the squared norm of $\frak p_m(x;a,b)$
(this is the standard notation, see \cite{Er, \S10.3}), and set
also $e_m=\frak p_m(1;a,b)$. It follows from  \cite{Er, \S10.3
(8)} and \tht{7.7} that
$$
a_m=\frac{k_{m+1}e_m}{k_m e_{m+1}}\,, \qquad
c_m=\frac{k_{m+1}k_{m-1}h_me_{m-1}}{k^2_mh_{m-1}e_{m+1}}\,,
$$
whence
$$
B(m,l)=\frac{k_{m+1}h_me^2_l}{k_mh_le_me_{m+1}}\,.
$$
Using the explicit values of the constants entering this formula
(see \cite{Er, \S10.8}) we get \tht{7.14}. \qed
\enddemo

\enddocument

\end